\documentclass[reqno]{amsart} 

\begin{document} 
\newtheorem{theorem}{Theorem}[section]
\newtheorem{lemma}{Lemma}[section]
\newtheorem{corollary}[theorem]{Corollary}
\newtheorem{example}[theorem]{Example}
\newtheorem{remark}{Remark}{}
\title[\hfil A study of asymptotic solutions]
{Asymptotic solutions of forced nonlinear second order differential equations and their extensions} 

\author[A. B. Mingarelli and K. Sadarangani \hfilneg]
{A. B. Mingarelli$^{(1)}$and K. Sadarangani$^{(2)}$ }  

\address{School of Mathematics and Statistics\\ 
Carleton University, Ottawa, Ontario, Canada, K1S\, 5B6\\ and \\ Departamento de Matem\'aticas, Universidad de Las Palmas de Gran Canaria, 
Campus de Tafira 
Baja, \\ 35017 Las Palmas de Gran Canaria, Spain.}
\email[A. B. Mingarelli]{amingare@math.carleton.ca, amingarelli@dma.ulpgc.es}
\address{Departamento de Matem\'aticas, Universidad de Las Palmas de Gran Canaria, 
Campus de Tafira 
Baja, \\ 35017 Las Palmas de Gran Canaria, Spain.}
\email[K. Sadarangani]{ksadaran@dma.ulpgc.es}

\date{}
\thanks{The first author is partially supported by a NSERC Canada Research Grant. Gratitude is also expressed to the Department of Mathematics of the University of Las Palmas in Gran Canaria for its hospitality during a visit there in Dec. 2004, and in 2006}

\begin{abstract}
Using a modified version of Schauder's fixed point theorem, measures of non-compactness and classical techniques, we provide new general results on the asymptotic behavior and the non-oscillation of second order scalar nonlinear differential equations on a half-axis. In addition, we extend the methods and present new similar results for integral equations and Volterra-Stieltjes integral equations, a framework whose benefits include the unification of second order difference and differential equations. In so doing, we enlarge the class of nonlinearities and in some cases remove the distinction between superlinear, sublinear, and linear differential equations that is normally found in the literature. An update of papers, past and present, in the theory of Volterra-Stieltjes integral equations is also presented.
\end{abstract}

\maketitle

\section{Introduction}

We present in this paper results pertaining to the nonlinear differential equation
\begin{equation}
\label{non0}
y''(x) + F(x, y(x)) = g(x), \quad x \in I = [x_0, \infty), x_0 \geq 0,
\end{equation}
where $F: \mathbb{R}^+ \times \mathbb{R} \longrightarrow \mathbb{R}$ is a general nonlinearity on which we will impose mostly criteria of integral type and $g(x)$ is given. Our main interest lies in the formulation of results regarding the non-oscillation and asymptotic behavior of its solutions. Some of the results will then be formulated for pure integral equations and ultimately for Volterra-Stieltjes integral equations (see \eqref{nlvs1}) and Volterra-Stieltjes integro-differential equations, that is, in the linear case, equations of the form
\begin{equation}\label{nlvs11}
y^{\prime}(x) = y^{\prime}(0) - \int_{0}^{x}y(t)\,d\sigma(t),
\end{equation}
and, in the nonlinear case, equations of the form
\begin{equation}\label{nlvs21}
y^{\prime}(x) = y^{\prime}(0) - \int_{0}^{x}F(t,y(t))\,d\sigma(t),
\end{equation}
where $\sigma$ is generally a function locally of bounded variation on $I$ and the resulting integrals are understood in the Riemann-Stieltjes sense. An advantage of the more general framework suggested by say, \eqref{nlvs11}, above is that one can incorporate corresponding theorems for three term linear recurrence relations such as
\begin{equation}
\label{31trr}c_ny_{n+1}+c_{n-1}y_{n-1} + b_ny_n = 0,\quad n\in \mathbb{N},
\end{equation}
and its nonlinear versions or, equivalently, second order linear difference equations such as
\begin{equation}
\label{2bdef}
{{\Delta}^2}{y_{n-1}}+b_ny_n = 0, \quad n\in \mathbb{N}.
\end{equation}
and its nonlinear analogs, as {\it corollaries} so that no new proof is  required to obtain the discrete analogs. \\

We recall that a solution of a real second order differential equation 
is said to be {\it oscillatory} on $[x_0, \infty)$ provided it exists on a semi-axis and it has arbitrarily large zeros on that semi-axis. If the equation has at least one non-trivial
solution with a finite number of zeros it is termed {\it non-oscillatory}. Recent work in asymptotics of \eqref{non0} has dealt primarily with pointwise criteria on both $F$ and $g$ sufficient for the asymptotic linearity of at least one solution (e.g., \cite{Con},\cite{mr},\cite{mr2},\cite{Yin},\cite{Zhao}, \ldots.)  On the other hand, {\it integral type} criteria cover by their very nature a wider collection of nonlinearities and we strive to obtain such criteria throughout. Thus, in Section~\ref{sec2} we give more general integral type criteria on $F$ which are sufficient for the existence of an uncoutable family of solutions of the unforced equation \eqref{non0}. This extends the validity of the results presented in Dub\'{e}-Mingarelli \cite{dm}. In addition, we note that our criteria of integral type such as \eqref{tkt} and \eqref{atk} below are extended over the whole half-line (that is we obtain {\it global existence}, see \cite{mr}) rather than local existence or existence for sufficiently large values of the variable. In this regard, see \cite{mps} for an extensive complete study of a specific nonlinear equation and \cite{mr} for a bibliographical study of unforced equations of the form $y^{\prime\prime}(x) + F(x, y(x), y^{\prime}(x)) = 0.$  For results which compare the non-oscillatory behavior of forced equations of the form \eqref{non0} with those of the associated unforced equation, \eqref{non} below, and possible equations with delays, we refer the reader to \cite{abu}, \cite{da} and \cite{ka}.\\

One should not forget that even though the literature is filled with sufficient criteria for oscillation/non-oscillation of unforced equations like
\begin{equation}
\label{non}
y''(x) + F(x, y(x)) = 0, \quad x \in I,
\end{equation}
in some cases, classical methods can actually be superior to the use of such fixed point theorems for the determination of the oscillatory character of an equation. For example, consider the equation $$y^{\prime\prime} + \frac{y\, \cos{2y^3}}{4(x+1)^2} =0, \quad x\in I,$$
whose nonlinearity fails to comply with the conditions of Nehari's theorem \cite{zn}, Atkinson's theorem \cite{1955}, the Coffman and Wong results in \cite{jsww3}, \cite{jsww4} and other more recent theorems. However, every solution of this equation is non-oscillatory as can be gathered by comparison with a non-oscillatory Euler equation (and use of Sturm's comparison theorem \cite{jcf}).\\

Next, we note that the use of maximum principles allows for an easy understanding of the oscillatory nature of an equation like \eqref{non}. For example, if in some interval $[a,b)$ (finite or not), we have $y\in C^2[a,b)$ and $y^{\prime\prime}(x) >0$ (or $y^{\prime\prime}(x) <0$) then $y(x)$ can have at most two zeros there. Thus, whenever a solution $y\in C^2$ of \eqref{non} satisfies $y^{\prime\prime}(x) \neq 0$, for $x\in I$, or, more generally, for all sufficiently large $x$, then we have non-oscillation on $I$. This explains the non-oscillatory character of equations like the Painlev\'{e} I, where $F(x,y) = -6y^2+x$, for $x > 0$ (see Hille~\cite{hille} or Ince~\cite{ince}). It follows that if $F$ is continuous and $F(x,y) < 0$ for all sufficiently large $x$ and all $y$, then we always non-oscillatory solutions on $I$ (and only such solutions). Thus, the only  interesting cases with regards to oscillations are those for which ultimately either $F(x,y) >0$ on its domain or $F(x,y)$ takes on both signs there. This motivates the main assumptions we will be making throughout.\\

As can be expected, introduction of the forcing term $g$ and its double primitive $f$, i.e.,  a function $f$ such that $g(x) = f^{\prime\prime}(x)$, can alter the original asymptotics. Loosely speaking, the case where $F$ dominates $g$ at infinity leads to solutions asymptotic to a double primitive of $g$ (see Section~\ref{sec3}). If $g$ is small in comparison to $F$, itself sufficiently small at infinity, then asymptotically linear solutions persist (see Section~\ref{ieq}). Motivated by Atkinson~\cite{fva} we introduce a novel necessary and sufficient condition for the existence of a solution of an integral equation of the form 
\begin{equation*}
y(x) = f(x) - \int_{x}^{\infty}(t-x)F(t,y(t))\,dt, \quad x\geq x_0
\end{equation*}
in terms of associated solutions of differential inequalities (Theorem~\ref{th4}). Ramifications of this result are noted and classical methods are used to obtain criteria for every solution of \eqref{non0} to be non-oscillatory. We also present an extension of Nehari's necessary and sufficient condition for non-oscillation [\cite{zn}, Theorem I], and Coffman and Wong's version \cite{jsww4} of the same in terms of solution asymptotics. We then proceed to a corresponding study of Volterra-Stieltjes integro-differential equations in Section~\ref{vse} and give conditions similar but more general than those in the previous sections. Finally, we apply this theory to obtain results for nonlinear three-term recurrence relations (or nonlinear second order difference equations). In addition, we give a long needed update of the theory of Volterra-Stieltjes integral equations in our Introduction to Section~\ref{vse}. For the purpose of clarity of exposition, we also proceed throughout the paper in order of increasing generality and leave the proofs until the very last section.

\section{Asymptotic results for nonlinear differential equations}

\subsection{Asymptotically linear solutions}
\label{sec2}
The present technique invokes a version of Schauder's fixed-point theorem
 and measures of non-compactness and is based, as in \cite{dm}, on the simple premise that  
in the variables separable case, the nonlinearity in the dependent variable
$y$ in \eqref{non} maps a given compact interval back into 
(and not necessarily onto) itself. For the rudiments of the notions of a measure of non-compactness and their applications, see the book by Bana\'{s} and Goebel \cite{bg}.

In the sequel, the space $ BC(\mathbb{R}^+)$ represents the space of all real bounded continuous functions defined on $\mathbb{R}^+$. For given $a \geq 0, b> 0$, we consider the space
\[
Y =\{u\in C[0,\infty): {\sup_{t \geq 0} \frac{|u(t)|}{at+b}} <\infty\}.
\]
Obviously, $Y$ is a vector space over $\mathbb{R}$. Now, for $u\in Y$ the quantity 
\[
\|u\|_Y =\sup_{t \geq 0} \frac{|u(t)|}{at+b}
\] is a norm on $Y$. Consideration of the mapping 
\[
\begin{array}{rlcl}
\Psi:& Y& \longrightarrow & BC(\mathbb{R}^+)\\
&u& \mapsto & \Psi (u)(t)= \displaystyle\frac{u(t)}{at+b} 
\end{array}
\]
shows that $\Psi$ is a linear operator and, moreover, $\Psi$ is an onto isometry. Consequently, as $BC(\mathbb{R}^+)$ is complete, $Y$ is 
a Banach space isometric to $BC(\mathbb{R}^+)$. More generally, for a given positive continuous function $p$, the space, $C_p$, of all tempered continuous functions (see [\cite{bg}, p.45]) consisting of all real-valued functions $ u \in C[0, \infty)$ such that $\sup_{t \geq 0} |u(t)|p(t) < \infty,$ is a Banach space. 

\begin{theorem} \label{th}
Let $a \geq 0, b>0$, and $X = \{ u \in Y |:  0 \leq u(t) \leq {at+b},\text{ for all }
 t \geq 0\}$. Assume that 
$F: \mathbb{R}^+ \times \mathbb{R}^+ \to \mathbb{R}^+$ is continuous and that for any 
$u \in X$, 
\begin{equation} \label{atk}
\int_{0}^{\infty} t\, F(t, u(t))\, dt \leq b.
\end{equation}
In addition, we assume that there exists a function $k:\mathbb{R}^+ \to \mathbb{R}^+$ with 
\begin{equation} \label{kt}
\int_{0}^{1}\, k(t)\, dt  < \infty,
\end{equation}
\begin{equation} \label{tktt}
\int_{0}^{1} t\, k(t)\, dt  < \infty,
\end{equation}
and
\begin{equation} \label{teek}
\int_{0}^{\infty} t^2\, k(t)\, dt  < \infty.
\end{equation}
and such that for any $u, v \in \mathbb{R}^+$,  
\begin{equation} \label{lip}
| F(t, u) - F(t, v) |  \leq k(t) | u-v|, \quad t \geq 0.
\end{equation}
Then \eqref{non} has a positive (and so non-oscillatory) asymptotically linear solution 
on $[0, \infty)$, i.e., $y(x) = ax+b + o(1),$ as $x \to \infty$.
\end{theorem}

\begin{remark}\label{uno} Note that \eqref{teek} does not necessarily imply neither \eqref{tktt} nor \eqref{kt}. However \eqref{kt}, \eqref{tktt}, and \eqref{teek} together do imply that 
$$\int_{0}^{\infty}\, t\, k(t)\, dt < \infty, \quad \int_{0}^{\infty}\, k(t)\, dt < \infty,$$ conditions that are used in various places in the proof.
\end{remark}

\begin{remark} We note in passing that if $a, b$ are chosen so that 
\begin{eqnarray}
\label{111}
\frac{1}{b}\, \max\{a,b\}\, \int_{0}^{\infty} t\,(t+1)\, k(t)\, dt < 1,
\end{eqnarray}
in the inequality \eqref{cont}, then $T$ is a contraction on $X$ and so the resulting fixed point is {\it unique} in $X$.
\end{remark}

\subsection{Asymptotic solutions in the forced nonlinear case}
\label{sec3}

In the sequel, the space $ BC([1,\infty))$ represents the space of all real bounded continuous functions defined on $[1,\infty)$. For a given function $g$ in \eqref{non} we assume that it has a second primitive $f : [1,\infty) \to \mathbb{R}$, such that for some $\delta > 0$
\begin{eqnarray}\label{ft} |f(x)| \geq \delta, \quad x \in [1,\infty)
\end{eqnarray}
a condition that we will return to and discuss at various points. Of course, since $f$ is continuous it is clear that \eqref{ft} implies that $f$ is of one sign on $[1,\infty)$, but the sign itself is of no concern to us here. Now, consider the vector space over $\mathbb{R}$ defined by
\begin{equation}
\label{spaceY}
Y =\{u\in C[1,\infty): {\sup_{x \geq 1} \frac{|u(x)|}{|f(x)|}} <\infty\}.
\end{equation}
Now, for $u\in Y$ the quantity 
\begin{equation}
\label{normY}
\|u\|_Y =\sup_{x \geq 1} \frac{|u(x)|}{|f(x)|}
\end{equation} 
is a norm on $Y$. Consideration of the mapping 
\[
\begin{array}{rlcl}
\Psi:& Y& \longrightarrow & BC([1,\infty))\\
&u& \mapsto & \Psi (u)(x)= \displaystyle\frac{u(x)}{|f(x)|} 
\end{array}
\]
shows that $\Psi$ is a linear operator and, moreover, $\Psi$ is an onto isometry. Consequently, as $BC([1,\infty))$ is complete, $Y$ is 
a Banach space isometric to $BC([1,\infty))$. More generally, for a given positive continuous function $p$, the space, $C_p$, of all tempered continuous functions (see [\cite{bg}, p.45]) consisting of all real-valued functions $ u \in C[x_0, \infty)$ such that $\sup_{x \geq x_0} |u(x)|p(x) < \infty,$ is a Banach space. \\

Let $F: [1, \infty) \times \mathbb{R} \to \mathbb{R}$ be continuous (not necessarily positive as in Section~\ref{sec2}) and assume that 
\begin{equation} \label{sf0}
\int_{1}^{\infty} s\, | F(s, 0) |\, ds < \infty.
\end{equation}
With $f$ as defined as in \eqref{ft}, we assume that there exists a function $k: [1,\infty) \to \mathbb{R}^+$ satisfying
\begin{equation} \label{sfk}
\int_{1}^{\infty}\,s\, |f(s)|\, k(s) \, ds  < \infty.
\end{equation}
An an additional restriction on both $F$ and $k$ we assume that for any $u, v \in \mathbb{R}$,  
\begin{equation} \label{lip01}
| F(x, u) - F(x, v) |  \leq k(x) | u-v|, \quad x \geq 1.
\end{equation}

Given such functions $F, k, f, g$ satisfying \eqref{ft}, \eqref{sf0}, \eqref{sfk} and \eqref{lip01} we consider the forced nonlinear equation \eqref{non} on the interval $I = [x_0, \infty)$ where $x_0$ is chosen so large that $x_0 \geq 1$ and for $x \geq x_0$,
\begin{equation} \label{t0}
 {\rm max} \left \{ \int_{x}^{\infty} \,(s-x)\, |f(s)|\, k(s) \, ds ,  \int_{x}^{\infty} \,(s-x)\, |F(s,0)| \, ds \right \} \leq \frac{\delta}{4},
\end{equation}
the finiteness of the integrals in \eqref{t0} being ensured on account of \eqref{sf0} and \eqref{sfk}. 

\begin{theorem} \label{th1}
Let the terms in \eqref{non} satisfy the conditions \eqref{ft}, \eqref{sf0}, \eqref{sfk}, \eqref{lip01} and \eqref{t0}. Consider \eqref{non} on $I = [x_0, \infty)$ where $x_0$ is defined as in \eqref{t0}.
Then \eqref{non} has a solution $y(x)$ satisfying
\begin{enumerate}
\item $y(x) \sim f(x)$ as $x \to \infty$ and actually, $y(x) = f(x) + o(1),$ as $x \to \infty$, and
\item 
\begin{equation}
\label{eqyf}
\sup_{x\in I} \frac{|y(x)|}{|f(x)|} \leq 2.
\end{equation}
\end{enumerate}
\end{theorem}

\subsection{Discussion} 

Since the proof of Theorem~\ref{th1} uses the contraction mapping principle, it follows that one can approximate the actual solution in question arbitrarily closely using a standard iterative technique.\\

The upper bound appearing in \eqref{eqyf} is by no means precise but will do for our purposes of obtaining global existence of solutions. Indeed, it is easily seen that one can modify the proof a little in order to find a non-uniform bound that depends  on $x_0$.\\

The unforced case $g(x) \equiv 0$ is included in our theorem and is reflected in the expression $f(x) =ax +b$ above, that is we obtain the existence of asymptotically linear solutions for \eqref{non0}. In this case Theorem~\ref{th1} extends the main results of Hallam \cite{tgh} for $n=1$, Dub\'{e}-Mingarelli \cite{dm}, and Mustafa-Rogovchenko \cite{mr}. It should be emphasized here that our conditions on the nonlinearity $F(x,y)$ and the forcing term $g(x)$ are essentially of {\it integral type} and not pointwise criteria as in most papers in the area, e.g., \cite{mr} is a recent one. In addition, Theorem~\ref{th1} provides an extension of some results in Atkinson \cite{fva} where, in addition, it is assumed that $F$ is positive and non-decreasing in its second variable (cf., also \cite{gp}), a condition we will return to occasionally.\\

Of course, since $f$ is continuous, \eqref{ft} implies that $f(x)$ is of {\it one sign} on the half-line $I$. Theorem~\ref{th1} then implies that the forced equation \eqref{non} is {\it non-oscillatory}. If \eqref{ft} is not satisfied then 
\begin{equation}
\label{liminff}
\liminf_{x\to \infty} |f(x)| = 0,
\end{equation}
a condition used often in many papers in conjunction with the questions under investigation here (\cite{fva}, \cite{ak}, \cite{ak2}, \ldots). In this respect, the condition \eqref{liminff} is known to furnish examples of oscillatory equations of the form \eqref{non}, cf., \cite{fva}, \cite{ak2}. In addition, necessary conditions for the existence of a positive solution of \eqref{non} under the assumption $f(x)> 0$, yet more restrictive conditions on the nonlinearity, may be found in [\cite{fva}, Section 4]. We note that \eqref{ft} and \eqref{sfk} together imply that 
\begin{equation}\label{sks}\int_{1}^{\infty}sk(s)\,ds < \infty,\end{equation}
 so that, as expected, one needs to ensure that the nonlinearity $F$ decreases quickly enough (see \eqref{lip01}) at infinity to ensure nonoscillation. Our results apply to linear problems with small forcing terms as well. The following example serves as illustration.

\begin{example} Let $F(x,y) = (1+y)/x^5$, $g(x)=1$, $x \geq 1$, in \eqref{non}. Choosing the double primitive $f(x) = x^2/2$, we see that $\delta =1/2$ is a suitable lower bound for $f(x)$ in \eqref{ft}. Note that \eqref{sf0}-\eqref{lip01} are all satisfied with the choice $k(x) = 1/x^5$. In addition, \eqref{t0} holds for all $x \geq x_0$ where $x_0=3$. Theorem~\ref{th1} now applies to show that the equation $$y^{\prime\prime}+(1+y)/x^5 = 1, \quad x\geq 3,$$ has a solution $y(x) \sim x^2/2$ as $x \to \infty$, defined by solving \eqref{tux} for its fixed point. This solution can actually be calculated using Bessel functions but its exact form is of no particular interest here. Successive approximations to it show that if we define $y_0(x)=1$, then $y_1(x) = x^2/2 - 1/{6x^3}$, and $$y_2(x) = \frac{x^2}{2} - \left (\frac{1}{4x} +\frac{1}{12x^3}-\frac{1}{252x^6}\right), \quad etc.,$$  the asymptotic nature of $y(x)$ can readily be ascertained.
\end{example}

A few more remarks on the case $g(x)=0$ in \eqref{non} are in order. Our condition \eqref{sf0} is compatible with Nehari's \cite{zn} necessary and sufficient condition for the existence of a bounded solution (albeit under additional assumptions on $F(x,y)$ such as positivity and monotonicity in its second variable). In this vein we can formulate the following immediate corollary for asymptotically constant solutions which does not assume neither the monotonicity nor the positivity of $F$.

\begin{corollary}Consider the equation \eqref{non} for $x \geq 1$. Let $F, k, \sigma $ satisfy \eqref{sf0}, \eqref{lip01}, \eqref{t0} and \eqref{sks}, for some $\delta=M>0$ and for all $x \geq x_0$. 
Then \eqref{non} has a solution satisfying $y(x) \to M $ as $x \to \infty$, and $|y(x)| \leq 2M$ for all $x \geq x_0$.
\end{corollary}

Similar additional results may be formulated for the case of asymptotically linear solutions and so are left to the reader. For an excellent survey up to the mid-seventies of nonlinear two term ordinary differential equations of Emden-Fowler type, see \cite{jsww}. \\

Next, we consider 
\begin{eqnarray}
\label{pert} y^{\prime\prime}(x) + F(x, y(x)) = g(x), \quad x\geq 0.
\end{eqnarray}
where $g(x) \equiv f^{\prime\prime}(x)$ and $f : \mathbb{R}^+ \to \mathbb{R}$ is not necessarily of one sign (as opposed to \eqref{ft}) but $f^{\prime\prime}\in C(\mathbb{R}^+)$. Next, we seek to find asymptotic theorems for equations of the form \eqref{pert} which may violate \eqref{ft}. The trade-off here is that we need that the nonlinearity be {\it positive}.

\begin{theorem}\label{th2} Let $f \in L^{\infty}( \mathbb{R}^+) \bigcap C^2( \mathbb{R}^+) $. Suppose that $F: \mathbb{R}^+ \times \mathbb{R}^+ \to \mathbb{R}^+$ is continuous and such that for some $b > 0$, \eqref{atk} is satisfied for any $u \in X$ where $X=\{ u \in BC(\mathbb{R}^+) : |u(x)| \leq \|f\|_{\infty} + b, x \geq 0\}$. Let $F$ satisfy 
\begin{equation} \label{lip}
| F(x, u) - F(x, v) |  \leq k(x) | u-v|, \quad x \geq 0.
\end{equation}
for any $u, v \in \mathbb{R}$, where $k:\mathbb{R}^+ \to \mathbb{R}^+$ is such that 
\begin{eqnarray}\label{tkt} \int_{0}^{\infty} t\, k(t)\, dt < 1.
\end{eqnarray}
Then \eqref{pert} has a solution $y(x)$ defined on $\mathbb{R}^+$ with $|y(x) - f(x)| \to 0$ as $x \to \infty$ and $\|y\|_{\infty} \leq \|f\|_{\infty} + b, x \geq 0.$
\end{theorem}

{\remark Uniqueness of the solution in Theorem~\ref{th2} may be lost in case we relax the requirement on $k$ as given by \eqref{tkt} to the integral being merely finite. In this case, a proof using measures of non-compactness such as the one in Theorem~\ref{th1} may be used to prove}

\begin{theorem}\label{th3} Let $f \in L^{\infty}( \mathbb{R}^+) \bigcap C^2( \mathbb{R}^+) $ and suppose that $F: \mathbb{R}^+ \times \mathbb{R}^+ \to \mathbb{R}^+$ is continuous and such that for some $b > 0$, \eqref{atk} is satisfied for any $u \in X$ where $X=\{ u \in BC(\mathbb{R}^+) : |u(x)| \leq \|f\|_{\infty} + b, x \geq 0\}$. Let $k:\mathbb{R}^+ \to \mathbb{R}^+$  satisfy \eqref{atk}, \eqref{lip} and 
\begin{eqnarray}\label{ttt} \int_{0}^{\infty} t\, k(t)\, dt < \infty.
\end{eqnarray}

Then \eqref{pert} has at least one solution $y(x)$ defined on $\mathbb{R}^+$ with $|y(x) - f(x)| \to 0$ as $x \to \infty$ and $\|y\|_{\infty} \leq \|f\|_{\infty} + b, x \geq 0.$
\end{theorem}
\vskip0.15in
\subsection{Discussion}

We make no claim as to {\it positivity} of the solution in question in either of Theorems~\ref{th2},\,\ref{th3}, since $f(x)$ may be of both signs, only that $y(x)=f(x) + o(1)$ as $x \to \infty$.  The following example illustrates this.

\begin{example} Consider the equation $$y^{\prime\prime} + F(x,y)= - \sin x, \quad x\geq 0,$$ with $F(x,y)= \lambda (x+1)^{-4}$ where $\lambda > 0$ is arbitrary but fixed and $b \geq \lambda/6$, where $b$ is defined in \eqref{atk}. For any constants $c_1, c_2$, we can choose the double primitive $f(x) = \sin x + c_1x+c_2$. The assumptions of Theorem~\ref{th2} are readily verified for our choice of $b$ and $F$. It follows that there is a solution of this equation such that $y(x) = f(x) + o(1)$ as $x \to \infty$. In fact the solution is given by $$y(x) = f(x) - \frac{\lambda}{6(x+1)^2},$$ for every $x$, from which the asymptotic estimate follows, as well as the {\it a priori} bound on the solution, namely, that $\|y\|_{\infty} \leq \|f\|_{\infty} + b$, valid for every $x \geq 0.$ Thus, choosing $c_1 =0, c_2=0$, we see that for large $x$ the solution will generally have both signs.
\end{example}

\section{Asymptotics for solutions of integral equations}
\label{ieq}

Motivated by Atkinson's paper \cite{fva} we produce a sharpening of the results in [\cite{fva}, Section 3] by studying {\it integral} inequalities. Our purpose is now to provide a formulations of some of the results of the previous sections to a wider framework, namely, that of integral equations, and ultimately to Volterra-Stieltjes integral equations with a view at obtaining discrete analogs for three-term recurrence relations.

Instead of beginning this study with a differential equation of the form \eqref{non0} we pass immediately to its integral equation counterpart, that is,
\begin{equation}
\label{ie1}
y(x) = f(x) - \int_{x}^{\infty}(t-x)F(t,y(t))\,dt, \quad x\geq x_0
\end{equation}
under various assumptions on the terms involved (after all, all our preceeding proofs were of this nature). Once again we strive to minimize the requirements on the forcing term, here, $f(x)$. In \cite{fva} this term is assumed to be small at infinity in the differential equation and differential inequality formulation. If $f\in C^2(I)$ we can recover results for the nonlinear equation \eqref{non0} by setting $g=f^{\prime\prime}$.\\

We will always assume that the ``forcing term" $f\in C(I)$ in \eqref{ie1} where $I$ as usual is of the form $I=[x_0,\infty)$, $x_0 \geq 0$, and uniformly bounded there. This last requirement will be denoted by the relation $f \in L^{\infty}(I)$, a minor abuse of notation. This is the only requirement we will impose upon $f$. The main result in this section follows:

\begin{theorem}
\label{th4} 
Let $f \in L^{\infty}(I)$ and suppose that the nonlinearity $F$ in $\eqref{ie1}$ satisfies
\begin{enumerate}
\item $F:I \times \mathbb{R}\to \mathbb{R}^+$ is continuous on this domain
\item $F(x, \cdot)$ is nondecreasing for every $x \in I$
\item For every $M > 0$, $$\int_{0}^{\infty}t\,F(t,M)\,dt < \infty$$
\item For every $y, z \in \mathbb{R}$ and every $x \in I$,
$$ | F(x,y)-F(x,z)| \leq k(x)|y-z|$$ where 
\item $$\int_{x_0}^{\infty}tk(t)\,dt < 1.$$
\end{enumerate}
Then \eqref{ie1} has a (continuous) solution $y \in L^{\infty}(I)$ if and only if there are two (continuous) functions $u, v \in L^{\infty}(I)$ such that $u(x) \leq v(x)$, $x \in I$,
\begin{equation}
\label{uf0}
u(x) \leq f(x) - \int_{x}^{\infty}(t-x)F(t,v(t))\,dt
\end{equation}
for $x\geq x_0$, and 
\begin{equation}
\label{vf0}
v(x) \geq f(x) - \int_{x}^{\infty}(t-x)F(t,u(t))\,dt
\end{equation}
for $x\geq x_0$.
\end{theorem}

{\remark \label{rem4} A corresponding result is valid for {\it positive} solutions $y$ of \eqref{ie1}. In this case $u(x) >0$ in Theorem~\ref{th4} and the positivity assumption on $F$ can be restated as $F(x,y) \geq 0$ for every $y \geq 0$, the remaining assumptions being the same. }

As a consequence we obtain a differential equations counterpart. Under the basic assumptions (1)-(5) of Theorem~\ref{th4}, we obtain

\begin{theorem} \label{cor01}Equation~\eqref{non0} with $g=f^{\prime\prime}$ has a solution $y$ with $y(x)=f(x) + o(1)$, $y^{\prime}(x)=f^{\prime}(x)+o(1)$, as $x \to \infty$ if and only if there exists two functions $u, v \in L^{\infty}(I)$ such that $u(x) \leq v(x)$, $x \in I$, satisfying \eqref{uf0} and \eqref{vf0},
\end{theorem}
\noindent{and }
\begin{theorem}\label{thx} Let $u, v, f$, $u(x) \leq v(x)$, for $x\in I$, be three twice continuously differentiable functions satisfying the differential inequalities $$u^{\prime\prime}(x) + F(x,v(x)) \leq g(x) \leq v^{\prime\prime}(x) + F(x,u(x)), \quad x \in I,$$ where $g=f^{\prime\prime}$. If, in addition, $u, u^{\prime}, v, v^{\prime}$ have vanishing limits at infinity, then \eqref{pert} has a positive solution $y \in L^{\infty}(I)$, with $y(x) \sim f(x)$ as $x \to \infty$. 
\end{theorem}

\begin{example}Consider the integral equation \eqref{ie1} with $F(x,y)=y/x^4$ and $f(x)=1+1/(6x^2)$ on $I=[1,\infty)$. Note that $y=1$ is a solution of \eqref{ie1} on $I$ and that for such $x$ all the conditions of  Theorem~\ref{th4} and Remark~\ref{rem4} are satisfied with the choice $k(x)=1/x^4$. The functions $u, v$, whose existence is guaranteed by this result, are given by $u(x)=1/2$ and $v(x) = 2$ for $x \in I$.
\end{example}
\vskip0.15in

\subsection{Discussion}

Our result gives, for a given (akin to `superlinear') nonlinearity, a necessary and sufficient condition for the existence of asymptotic solutions of a given type in terms of solutions to specific integral inequalities (according to \cite{jsww3} a superlinear $F$ is one which satisfies condition `2' in Theorem~\ref{th4}). The stated theorem gives uniqueness as well as uniform bounds for the required solution. It is likely that hypotheses `4' and `5' in the theorem may be relaxed albeit at the possible loss of uniqueness. Theorem~\ref{th4} appears to be new even when considered from the viewpoint of  second order nonlinear differential equations (as in, e.g., Theorem~\ref{cor01}).

As mentioned earlier, Nehari \cite{zn} gives a necessary and sufficient condition for the existence of a non-oscillatory solution of an equation akin to \eqref{non0} in terms of the nonlinearity. Under a {\it strong superlinear} condition that is, for some $\varepsilon > 0$ there holds
\begin{equation}
\label{sts}
y_2^{-\varepsilon}F(x,y_2) >y_1^{-\varepsilon}F(x,y_1) ,
\end{equation}
for $0< y_1<y_2<\infty$, he shows in \cite{zn} that if the nonlinear equation has the special form
\begin{equation}
\label{zn11} y^{\prime\prime}+yF(x, y^2) = 0, \quad x\geq 0,
\end{equation}
then it has a bounded non-oscillatory solution if and only if, for some $M >0$,
\begin{equation}
\label{zn2} \int_{}^{\infty}tF(t,M)\, dt < \infty,
\end{equation}
i.e., condition `3' in Theorem~\ref{th4} holds for {\it some} $M>0$ (and $x=0$).

His result was extended later by Wong \cite{jsww2} by relaxing the monotonicity condition on $F$ somewhat and taken up again by Coffman-Wong \cite{jsww3}, \cite{jsww4},  where further developments in a sublinear case were given (in particular, see the Table on p.123 in \cite{jsww2} for a useful visual display of known necessary and sufficient conditions for non-oscillation). A more precise version of Nehari's result can be found in the strong superlinear case in \cite{jsww4}, that is, \eqref{zn11} has a bounded {\it asymptotically linear} solution (the special case $f(x)=ax+b$ in our set up) if and only if \eqref{zn2} holds for some $M >0$.

Although we require that \eqref{ie1} be `superlinear' (i.e., condition `2' in Theorem~\ref{th4}) it need not be strongly so. On the other hand, we require that \eqref{zn2} hold for every $M >0$, but then we are also strengthening the conclusion. Indeed, our result is also valid for a wider class of equations, not only second order nonlinear differential equations. 

In a recent paper \cite{mr2}, the authors implicitly assume an integrability condition on $f$ and that this double primitive $f(x) \to 0$ as $x \to \infty$, in the spirit of \cite{fva} and \cite{ak}. The nature of such a decay condition on the forcing term is that the basic tenet  underlying the asymptotic behavior of a given nonlinear differential equation \eqref{non0} appears to be the interplay between the rate of decay of the nonlinearity as opposed to the rate of decay of the forcing term. For example, if the nonlinearity is small in the sense of the applicability of conditions (3-5) in Theorem~\ref{th4} and the forcing term is not ( e.g., perhaps not integrable on I) then solutions may be expected to be asymptotic to a double primitive of the forcing term. On the other hand, if both forcing term and nonlinearity are ``small" in some suitable sense then the solutions of \eqref{non0} may be expected to be asymptotically linear (or asymptotic to the solutions of the same equation with nonlinearity and forcing term omitted). This philosophy may be used in our basic understanding of nonlinear equation asymptotics on a half-line on account of the following unpublished result, reproduced here for completeness.

\begin{theorem} \label{atm}(Atkinson-Mingarelli, 1976, unpublished)
\end{theorem}Let $g:I\to (0,\infty)$ be continuous and satisfy
\begin{equation}
\label{tgt}\int_{x_0}^{\infty}t^i\, g(t)\, dt < \infty
\end{equation}
for $i=0,1$. We assume that $G:I \times \mathbb{R}\to \mathbb{R}^+$ verifies assumptions (1-2) in Theorem~\ref{th4} (i.e., $G$ is continuous on this domain and $G(x, \cdot)$ is nondecreasing for every $x \in I$). In addition, let $G_{x}(x,y) \equiv \partial{G}/ \partial{x}$ exist, be continuous and non-positive on the same domain. If
\begin{equation}
\label{tgKt} \int_{x_0}^{\infty}t\, G(t, Kt)\, dt < \infty
\end{equation}
for every $K > 0$, then every solution $y$ of
\begin{equation}
\label{non2} y^{\prime\prime} + yG(x,y) = g(x)
\end{equation}
is either asymptotically linear or $ y(x) < 0$, $y^{\prime}(x) \leq 0$, $y^{\prime\prime} \geq 0$ for all sufficiently large $x$. Either way, every solution is nonoscillatory.

{\remark The double primitive of $g$ does not enter the picture here as in Theorem~\ref{th4} since $g$ is already itself small, as evidenced by \eqref{tgt}.  By a solution $y$ that is ``asymptotically linear" we mean that the solution has the property that $y(x)= Ax+B+o(1)$ as $x \to \infty$, for some constants $A, B$, where the cases $A=0, B> 0$, $A > 0, B=0$, $A=B=0$ can all occur. This result is in the same spirit as Theorem~\ref{th4} above except for clear differences in the behavior of the nonlinearities involved. Despite these differences, these nonlinearities are each small at infinity thereby leading to the stated {\it linear} asymptotics. This theorem extends a theorem of [Nehari \cite{zn}, Theorem III]
}

The following complementary result to Theorem~\ref{atm} is for the case where $g$ is not integrable at infinity. In this case the nonlinearity is still small in comparison to the growth of $g$ in \eqref{fr} but now the forcing term is integrably large, so the solutions are asymptotic to a double primitive of the forcing term.

\begin{theorem} \label{atm2}(Atkinson-Mingarelli, 1976, unpublished)
\end{theorem}Let $g:I\to (0,\infty)$ be continuous and satisfy
\begin{equation}
\label{tgt2}\int_{x_0}^{\infty} g(t)\, dt = + \infty
\end{equation}
for $i=0,1$. We assume that $G:I \times \mathbb{R}\to \mathbb{R}^+$ verifies assumptions (1-2) in Theorem~\ref{th4} (i.e., $G$ is continuous on this domain and $G(x, \cdot)$ is nondecreasing for every $x \in I$). In addition, let $G_{x}(x,y) \equiv \partial{G}/ \partial{x}$ exist, be continuous and non-positive on the same domain. If for some double primitive $f$ (i.e., $g = f^{\prime\prime}$) and some $\varepsilon > 0$
\begin{equation}
\label{fr} \int_{x_0}^{x}\sup_{|u| \leq (1+\varepsilon)f} |F(t,u)|\, dt = o\left \{\int_{x_0}^{x}g(t)\, dt\right\},
\end{equation}
as $x \to \infty$, then every solution $y$ of
\begin{equation}
\label{non2} y^{\prime\prime} + yG(x,y) = g(x)
\end{equation}
is asymptotic to $f(x)$ as $x \to \infty$.

\subsection{Discussion}
\label{dis4}

One of the advantages in using classical methods over fixed point theorems is exhibited in Theorems~\ref{atm} and~\ref{atm2} above. For example, it is difficult to obtain {\it a priori} bounds such as \eqref{ener3} on the derivative or \eqref{yKx} on the solution $y(x)$ using fixed point theorems. Both techniques can be used interchangeably, preference being only a function of the conclusion desired, and on the nature of the hypotheses, nothing more.

Note, however, the absence of conditions such as \eqref{lip} or \eqref{tkt} in these two results, hypotheses that were deemed necessary in the proofs of most of the results in this section. Conditions such as \eqref{lip} and \eqref{tkt} could also be interpreted as being conditions on the rate of growth of $\partial{G}/\partial{y}$ in the domain under consideration.

\section{Asymptotic theory of nonlinear Volterra-Stieltjes integral equations}
\label{vse}
We begin this section with a non-exhaustive review and update of the results over the past 20 years in this fascinating area which can be used to unify discrete and continuous phenomena. The unification allows for the simultaneous study of both differential and difference equations of the second order and even includes equations that are, in some sense, {\it between} these two as we gather from the discussion that follows (and from the references). A  {\it Volterra-Stieltjes integral equation} is basically a Volterra integral operator on a space $X$ of suitable functions in which the integral appearing therein is a Stieltjes integral (in whatever sense it can be defined, more on this below). The prototype (linear) Volterra-Stieltjes integral equation that we use in this work is of the form $I \equiv 0 \leq x < b \leq\infty$,
\begin{equation}\label{vs1}
y(x) = y(0)+xy^{\prime}(0) - \int_{0}^{x}(x-t)y(t)d\sigma(t)\quad x \in I,
\end{equation}
where $\sigma: I \to \mathbb{R}$ is a function that is locally of bounded variation on $I$. A {\it solution} of \eqref{vs1} is an absolutely continuous function with a right-derivative that exists for each $x \in I$ and is locally of bounded variation on $I$. In this case, the integral in \eqref{vs1} can be understood in the {\bf Riemann-Stieltjes} sense and we take this for granted throughout this section. It is known that this formulation, which includes the use of a simple Riemann-Stieltjes integral, is adequate (and sufficient) for the unification purposes referred to above (see e.g., \cite{abmb} among other possible references).  Other frameworks that can be used as a unification tool for discrete and continuous nonlinear equations include the theory of {\it time scales}. However, we do not entertain these studies here (unless results overlap) and so refer the interested reader to e.g., \cite{mbo},\cite{lhe},\cite{erk},\cite{billur},\cite{shs}, and the references therein for further information. \\

Associated with \eqref{vs1} is the Volterra-Stieltjes integro-differential equation
\begin{equation}\label{vs2}
y^{\prime}(x) = y^{\prime}(0) - \int_{0}^{x}y(t)d\sigma(t)
\end{equation}
obtained by differentiating the equation \eqref{vs1}. The derivative appearing in \eqref{vs2} is now understood generally as a right-derivative and this function is locally of bounded variation on $I$. Local existence  and uniqueness of solutions of initial value problems associated with either \eqref{vs1} or \eqref{vs2}  and the basic theory of such equations, as we define them, was developed by Atkinson \cite{fvab}, and also continued by Mingarelli \cite{abmb}, Mingarelli and Halvorsen \cite{mhb} among others. The reader may also wish to consult the monographs of Corduneanu \cite{cc}, H\"{o}nig \cite{csh}, Schwabik {\it et al.} \cite{ssc3} and \cite{ssc1} for different approaches and generalizations of both the method and the context. We also refer the reader to Groh \cite{jg0} for an extensive list of references to this subject, some not included here. \\

Regarding the possible different interpretations of the nomenclature ``Volterra-Stieltjes integral equation" in the literature, the main ideas and developments depend strongly on the definition of the particular Stieltjes integral being used. Once this is in place, one can define a solution, develop a basic theory (ask about existence and uniqueness of solutions, continuous dependence on initial conditions, etc) and suggest additional applications. Thus, studies which have depended upon the use of the {\bf Kurzweil-Henstock integral} in \eqref{vs1} and equations like it include, and are not restricted to, Kurzweil \cite{jk}, Schwabik \cite{ssc4}, Tvrd\'y \cite{mtu1}, Federson-Bianconi \cite{feb}. One of the very first  researchers in this area, Martin~\cite{rhm} uses the {\bf Cauchy right-integral} while the {\bf  Dushnik integral} appears in  the papers by H\"{o}nig~\cite{csh}, \cite{csh3} and Dzhgarkava~\cite{ddt2}. The Stieltjes integral, when viewed as a {\bf Lebesgue-Stieltjes} integral, makes its appearance in Ding-Wang \cite{xqd1}, \cite{xqd2} and Caizhong~\cite{lcz}. Finally, but not exhaustively, Dressel~\cite{fgd} uses the {\bf Young} integral. There may be overlap between some of these definitions as they have developed in the past century, but the aim is to show that different integral definitions may produce different applications and may account for the large literature on this subject.\\ 

Various {\bf abstract formulations} of such equations can be found in the works of Helton~\cite{bwh} who considered such equations over rings, Ashordia~\cite{mash}, \cite{mash1}, Dzhgarkava~\cite{ddt}, H\"{o}nig~\cite{csh4},\cite{csh5}, Ryu~\cite{ksr}, Schwabik~\cite{ssc0}, Travis~\cite{cct} and Young~\cite{dfy}.
{\bf Controllability} of equations of this general type has been studied by Barbanti~\cite{lb}, Dzhgarkava \cite{ddt2}, Groh~\cite{jg2}, Yong~\cite{jmy}, and Young~\cite{dfy2}. Investigations related to {\bf systems of equations} include those of Gopalsamy {\it et al.} \cite{kgo}, Herod~\cite{jvh}, Hildebrandt~\cite{thh}, Hinton~\cite{dbh}, H\"{o}nig~\cite{csh2}, Schwabik~\cite{ssc},  and Wheeler~\cite{rlw}. For the relationship between Volterra-Stieltjes integral equations and their applications to {\bf difference equations}  (or recurrence relations) see Atkinson~\cite{fvab},  Mingarelli~\cite{abmb}, Mingarelli-Halvorsen \cite{mhb}, Petrovanu~\cite{dp}, and Schwabik~\cite{ssc4}. \\

The study of {\bf scalar} Volterra-Stieltjes integral equations and subsequent qualitative, quantitative, and spectral theory can be found in the works by Banas {\it et al}~\cite{bg2}~\cite{bg3}~\cite{bg4}~\cite{bg5}~\cite{bg6}~\cite{bg7}~\cite{bg8}, Caballero {\it et al}~\cite{crs}, Cao~\cite{zjc}, Cerone-Dragomir~\cite{pce}, Chen~\cite{hyc}~\cite{hyc1}, El-Sayed~\cite{wes}, Gibson~\cite{wlg}, Gil' and Kloeden~\cite{mig}~\cite{mig2}, Hu~\cite{qyh}, Jiang~\cite{zmj}~\cite{zmj0}, Lou~\cite{lou}, Marrah and Proctor~\cite{marr}, Mingarelli~\cite{abm0}~\cite{abmb}, Mingarelli and Halvorsen~\cite{mhb}, Parhi~\cite{np}, Randels~\cite{wcr}, Schwabik~\cite{ssc2}, Spigler and Vianello \cite{rsv1}~\cite{rsv}, Tritjinsky~\cite{tri}, Wang~\cite{zwa}, Wong and Yeh~\cite{wy0}~\cite{wy}.\\

Finally, there is a relationship which has seen little follow-through in the past 70 years or so since its beginnings. Basically, one asks about a relationship between the notion of a generalized derivative ({\it \`{a} la Feller})  of the form $$ - \frac{dy^{\prime}}{d\sigma} = f(x), \quad x\in [0,b],$$ and the integro-differential equation (cf., \eqref{vs2} above), $$y^{\prime}(x) = c - \int_{0}^{x}f(t)\,d\sigma(t),$$ where $c$ is a constant and $\sigma$ is an non-decreasing (actually increasing) function defined on $[0,b]$ with an appropriate Stieltjes integral. This approach was pioneered by the probabilist W. Feller~\cite{wfe} (see the references in [\cite{abmb}, p.316]) and the resulting theory, found in many papers in probability (not all quoted here) now includes the keywords: Feller derivatives, Krein-Feller operators, Generalized differential operators, etc., see Jiang~\cite{zmj2}, Albeverio-Nizhnik~\cite{sal}, Fleige~\cite{anf}, Mingarelli~\cite{abmb}.\\

That there is an equivalence between these last two displays should not be surprising yet the lines of development of the resulting theories seem to have diverged over the years with each equation taking on a life of its own, so to speak. For papers dealing with generalized differential expressions see Groh~\cite{jg1}, Jiang~\cite{zmj2}, Volkmer~\cite{hv}, and Mingarelli~\cite{abmb} among others. We emphasize here the importance of the contributions of I.S. Ka\v{c}, M.G. Kre\v{i}n and H.K. Langer to the study of the spectral theory of the operators associated with the generalized differential expressions above. Although these references are not included here specifically for reasons of length, we refer the interested reader to the more than 80 historical references in Mingarelli~\cite{abmb}, in addition to those in Atkinson~[\cite{fvab}, pp.529-533], the list of references in Fleige~\cite{anf} and the references contained within each of the articles mentioned in the bibliography below. Altogether these should give the reader an essentially complete view of this vast field as of today.\\

\subsection{Asymptotically linear solutions of nonlinear equations}
\label{alis}
The asymptotic theory of solutions of equations of the form \eqref{vs1} or \eqref{vs2} is still in its infancy with few basic results in existence in the literature. In the linear case \eqref{vs1} we can cite Theorem~12.5.2 in Atkinson~\cite{fvab}. Fewer are specific results dealing with the nonlinear case \eqref{nlvs1} or \eqref{nlvs2}. One such result may be found in [\cite{abmb}, Theorem 2.3.1] in the case where $F(x,y) = p(x)q(y)$, a result which extends Butler's necessary and sufficient condition for non-oscillation \cite{Bu}. In the remaining sections we produce extensions of the results in the previous sections to this framework along with some possible refinements.

\begin{equation}\label{nlvs1}
y(x) = y(0)+xy^{\prime}(0) - \int_{0}^{x}(x-t)F(t,y(t))\,d\sigma(t)
\end{equation}

\begin{equation}\label{nlvs2}
y^{\prime}(x) = y^{\prime}(0) - \int_{0}^{x}F(t,y(t))\,d\sigma(t)
\end{equation}
Using the methods in Atkinson~[\cite{fvab}, Chapter 12], one can readily prove the existence and uniqueness of solutions of initial value problems for equations of the form \eqref{nlvs1} or \eqref{nlvs2} under a locally Lipschitz condition on the continuous nonlinearity $F$. Recall that a solution of \eqref{nlvs1} (resp.\eqref{nlvs2}) is an absolutely continuous function such that its right derivative exists at every point of $I$ and $y(x)$ (resp.$y^{\prime}(x)$) satisfies the equation \eqref{nlvs1} (resp.\eqref{nlvs2}) at every point in $I$. Unless otherwise specified we always assume the minimum requirement that such solutions exist and are unique. In some cases below we actually get existence, uniqueness and asymptotic limits as a by-product of the techniques used.\\

For simplicity of notation we will assume hereafter, unless otherwise specified, that the interval $I$ in question is $I=[0,\infty)$, but it could well be any half-line, of the form $I=[x_0, \infty)$, with minor changes throughout (obtained by a change of independent variable). Our first general result is a counterpart of Theorem~\ref{atm} for asymptotically linear solutions of \eqref{nlvs1}.

\begin{theorem}
\label{th5}
Let $\sigma : I \to \mathbb{R}$ be right continuous and locally of bounded variation on $I$. Suppose that the nonlinearity $F$ in $\eqref{nlvs1}$ satisfies
\begin{enumerate}
\item $F:I \times \mathbb{R}\to \mathbb{R}^+$ is continuous on this domain
\item $F(x, \cdot)$ is nondecreasing for every $x \in I$
\item For some $M > 1$, $$\int_{0}^{\infty}F(t,Mt)\,|d\sigma(t)| < \infty$$
\end{enumerate}
Then \eqref{nlvs1} has an asymptotically linear solution, viz., a solution $y$ with $y(x) = Ax+B +o(1)$ as $x \to \infty$ for some appropriate choice of real numbers $A, B.$
\end{theorem}

\begin{remark} The assumption that $F$ is nondecreasing in its second variable may be weakened at the expense of additional smoothness as a function of that variable (e.g., a Lipschitz condition of type \eqref{lip} and \eqref{ttt} as we have seen above) and use of a fixed point theorem as the next result shows.
\end{remark}

\begin{theorem}
\label{th6}
Let $\sigma:I\longrightarrow\mathbb{R}$ be a
non-decreasing right-continuous function, $F: I \times\mathbb{R}^{+}\longrightarrow\mathbb{R}^{+}$ be continuous such that for some $M > 0$, 
\begin{itemize}
\item[$(a)$]$\displaystyle\int_{0}^{\infty}t\,F(t,y(t))\,d\sigma(t)\leq M,$ \quad{for $y\in X$,} 
\end{itemize}
where $X=\{y\in \mathcal{C}(I): 0\leq y(x)\leq M,\,\,\,\,\, x\in I\}$,
\begin{itemize}
\item[$(b)$]$|F(x,u)-F(x,v)|\leq k(x)|u-v|$, \quad{$x \in I$, $u, v \in 
\mathbb{R^+}$}
\end{itemize}
where $k: I \to \mathbb{R}^{+}$ is continuous and
\begin{itemize}
\item[$(c)$]$\displaystyle\int_{0}^{\infty}t\,k(t)\,d\sigma(t)<\infty$.
\end{itemize}
Then the Volterra-Stieltjes integro-differential equation \eqref{nlvs2} has a monotone increasing solution $y(x)$ with $0 \leq y(x) \leq M$ for $x \in I$ and $y(x) \to M$ as $x \to \infty$.
\end{theorem}

It appears at first sight as if condition (a) in Theorem~\ref{th6} may be difficult to verify. However, the following simple corollary shows that pointwise estimates on $F(x,y)$ can be used to imply the same conclusion.

\begin{corollary} \label{coro} 
Assume that $F, \sigma$ are as in Theorem~\ref{th6}. Let $M > 0$ and let 
\begin{equation}\label{atko}
 F(x, y)  \leq  p(x) q(y),\quad x \geq 0, y \in \mathbb{R}^+,
\end{equation}
for some function $q$, where $q : [0, M] \to [0, M] $ is continuous on 
$[0,M]$. Let $p\in C[0,\infty)$ and suppose that
\begin{equation}
\label{effo}
 \int_{0}^{\infty}t\, p(t)\, d\sigma(t)\, \leq 1.
\end{equation} 
Assume further that there exists a function 
$k:\mathbb{R}^+ \to \mathbb{R}^+$ such that $k$ is continuous and 
\begin{equation*}
 \int_{0}^{\infty}t\, k(t)\, d\sigma(t)\, < 1
\end{equation*}
such that for any $u, v \in \mathbb{R}^+$, we also have 
\begin{equation*} 
| F(x, u) - F(x, v) | \leq  k(x) | u-v|, \quad x \geq 0\,.
\end{equation*}
Then \eqref{nlvs2} has a positive (and so non-oscillatory) monotone solution 
on $I$ such that $y(x) \to M$ as $x \to \infty$.
\end{corollary}

\subsection{Discussion}

Note that if $\displaystyle\int_{0}^{\infty}t\,F(t,0)\,d\sigma(t)<\infty$ then this condition, along with assumptions (b) and (c) in the theorem together imply (a). In particular, (a) is satisfied if $F(x,0)=0$ for every $x \in I$.

As in the differential equation case before, if $\displaystyle\int_{0}^{\infty}t\,k(t)\,d\sigma(t)<1$ then relation \eqref{ayxaz} in its proof gives us
$$
\|Ay-Az\|_{\infty}\leq\|x-y\|_{\infty}\int_{0}^{\infty}t\,k(t)\,d\sigma(t),
$$
and the Banach contraction mapping theorem applies immediately to gives us existence and uniqueness of the solution of our integro-differential equation. Note the similarity between hypothesis (3) in Theorem~\ref{th5} and assumption (a) in Theorem~\ref{th6}: In condition (a) the integrand involves a class of functions all bounded by the constant $M$, whereas in hypothesis (3) the ``class of functions" is replaced by the class of linear functions of the form $Mt$. In the former case there are asymptotically constant solutions while, in the latter case, there are asymptotically linear solutions. This is reflected in the form of the respective assumptions. Indeed, since the constant function $y(x)=M$ is in $X$, assumption (a) {\it includes} the condition 
\begin{equation}\label{tftm}
\int_{0}^{\infty}t\,F(t,M)\, d\sigma(t) \leq M.
\end{equation}
Next, Theorem~\ref{th6} gives the existence of an asymptotically constant solution whenever there exists a constant $M >0$ satisfying condition (a) (the other two assumptions being independent of $M$ we assume as implicitly verified). Thus, if (a) is assumed for {\it each} $M > 0$ then there is an asymptotically constant solution tending to that limit, $M$. A similar observation appplies to Theorem~\ref{th5}. A moment's reflection shows that if, in addition, we assume that $F(x, \cdot)$ is non-decreasing for each $x \in I$, then the existence of a solution can be obtained satisfying the improved estimate $$ M - \int_{0}^{\infty}t\,F(t,M)\, d\sigma(t) \leq y(x) \leq M$$ in Theorem~\ref{th6}. Since \eqref{tftm} holds for that $M$, the left hand side is non-negative. Since \eqref{tftm} is reminiscent of Nehari's criterion \cite{zn} for the existence of bounded nonoscillatory solutions, it is of interest to investigate the validity of this criterion in this more general setting and this is the subject of the next result.

\begin{lemma}
\label{lem1}
Let $\sigma:I\longrightarrow\mathbb{R}$ be a
non-decreasing right-continuous function, $F: I \times\mathbb{R}^{+}\longrightarrow\mathbb{R}^{+}$ be continuous and such that for some $M > 0$, 
\begin{equation}
\label{zn1}
\int_{0}^{\infty}t\,F(t,M)\,d\sigma(t)< \infty. 
\end{equation}
Then every eventually positive solution of the Volterra-Stieltjes integro differential equation \eqref{nlvs2} is either of the form $y(x)\sim Ax$ as $x \to \infty$  for some constant $A\neq 0$ or $y(x)/x \to 0$ as $x \to \infty$. 
\end{lemma}

We now formulate an analog of Nehari's necessary and sufficient criterion \cite{zn} for the existence of a bounded nonoscillatory solution of our equation (recall that a solution $y$ of \eqref{nlvs1} or \eqref{nlvs2} is said to be {\it nonoscillatory} provided $y(x) \neq 0$ for all sufficiently large $x$).

\begin{theorem}
\label{th7}
Let $\sigma:I\longrightarrow\mathbb{R}$ be a
non-decreasing right-continuous function, $F: I \times\mathbb{R}^{+}\longrightarrow\mathbb{R}^{+}$ be continuous and non-decreasing in its second variable (i.e., $F(x, y)$ is nondecreasing in $y$ for $y>0$, for each $x\in I$). Then \eqref{nlvs2} has bounded eventually positive solutions if and only if \eqref{zn1} holds for some $M>0$.
\end{theorem}

\begin{corollary}\label{cor001} Let $\sigma$ be as in Theorem~\ref{th7}, $G: I\times \mathbb{R}^+\to \mathbb{R^+}$ be continuous and positive in $I\times \mathbb{R}^+$. In addition, let $G(x,y)$ be nondecreasing for every $y>0$, $x\in I$. Then 
\begin{equation}\label{nlvs2p}
y^{\prime}(x) = y^{\prime}(0) - \int_{0}^{x}y(t)G(t,y^2(t))\,d\sigma(t)
\end{equation}
has bounded nonoscillatory solutions if and only if there holds 
\begin{equation}
\label{zn18}
\int_{0}^{\infty}t\,G(t,c)\, d\sigma(t) < \infty,
\end{equation}
for some $c>0$.
\end{corollary}

\noindent{T}he proof of the next result is an immediate consequence of the theorem.

\begin{corollary}\label{cor02} Let $\sigma, F$ be as in Theorem~\ref{th7}. Then \eqref{nlvs1} has asymptotically constant positive solutions if and only if \eqref{zn1} holds for some $M>0$.
\end{corollary}

Of course, Corollary~\ref{cor02} deals with {\it bounded} solutions of \eqref{nlvs1}. An analogous result for possibly unbounded solutions follows (although strong superlinearity \eqref{sts} is to be imposed).

\begin{theorem}\label{th8} Let $\sigma$ be as in Theorem~\ref{th7}. Assume that $F: I\times \mathbb{R}^+\to \mathbb{R^+}$ is continuous and positive in $I\times \mathbb{R}^+$. In addition, let $F$ satisfy the strong superlinearity condition \eqref{sts} for $\varepsilon =0$, as well as for some $\varepsilon > 1$. Then \eqref{nlvs1} has an eventually positive solution if and only if  \eqref{zn1} holds for some $M>0$.
\end{theorem}

\subsection{Discussion} Theorem~\ref{th7} is an improvement of Nehari's theorem \cite{zn} to the framework of Volterra-Stieltjes integral equations \eqref{nlvs1}, or Volterra-Stieltjes integro-differential equations \eqref{nlvs2}.  Although Nehari's theorem [\cite{zn}, Theorem I] was stated for equations of the form \eqref{zn11}, we choose the more general form stated here, with an arbitrary nonlinearity (this explains the apparently odd restriction on $\varepsilon > 1$ rather than $\varepsilon > 0$ as in the original Nehari result). As pointed out in the proof of Corollary~\ref{cor001} the form \eqref{zn11} is actually guided by the wish that both $y$ and $-y$ be solutions of the same equation. Nehari's theorem as such is actually a special case of Corollary~\ref{cor001} with $\sigma(t)=t$ throughout. The integral equation \eqref{nlvs2p} then produces a differential equation of the form \eqref{zn11} (since the indefinite integral is continuously differentiable). Indeed, Corollary~\ref{cor02} (via the techniques in the proof Corollary~\ref{cor001}) also includes an extension of Nehari's theorem by Coffman and Wong \cite{jsww4}, [\cite{jsww3}, Theorem E]. \\

\noindent{T}he Volterra-Stieltjes framework provides for recurrence relation (discrete) analogs or even intermediate mixed type integro differential equations as a direct consequence (see the next Section for applications). In addition, Corollary~\ref{cor02} shows that the sufficiency of the proof of Theorem~\ref{th7} actually provides a criterion for the existence of asymptotically constant solutions of either \eqref{nlvs1} or \eqref{nlvs2}. As we gather from the proof of said theorem, we can choose the asymptotic limit $A$ appearing in \eqref{ii0} to be any number between $(0, M)$, where the $M$ appears in \eqref{zn1}. It follows that if \eqref{zn1} is valid for every $M>0$ then \eqref{nlvs1} has solutions whose limits can be any prescribed positive number.\\

\noindent{T}heorem~\ref{th8} includes a slight modification of an additional result of Coffman and Wong [\cite{jsww4}, Section 6]. Observe that, if the solution in the necessity of Theorem~\ref{th8} is unbounded, then \eqref{zn1} must hold for {\it every} $M>0$, just as in the case of ordinary differential equations, cf., \cite{jsww4}. That is, the existence of at least one unbounded eventually positive solution of \eqref{nlvs1} implies the convergence of the integral \eqref{zn1}, not only for the $M$ in question, but for every $M>0$ (see also Lemma~\ref{lem2} below in this regard). \\ 

In order not to restrict ourselves only to the study of asymptotically constant solutions of either \eqref{nlvs1} or \eqref{nlvs2}, we now present further results relating to asymptotically linear solutions. Lemma~\ref{lem2} below complements Theorem~\ref{th5} above.\\

\begin{lemma}
\label{lem2}
Let $\sigma$ be right-continuous and nondecreasing on $I$, $F: I\times \mathbb{R}^+\to \mathbb{R^+}$ be continuous and positive in $I\times \mathbb{R}^+$. In addition, let $F(x,y)$ be nondecreasing for every $y>0$, $x\in I$. If either \eqref{nlvs1} or \eqref{nlvs2} has a solution $y(x) \sim Ax+B$ as $x \to \infty$, where $A>0$, $B$ are constants, then 
\begin{equation}
\label{tmt}
\int_{0}^{\infty} F(t, Mt)\, d\sigma(t) < \infty
\end{equation}
for some $M > 0$.
\end{lemma}

\begin{remark}
Incidentally, this proof also shows that the existence of at least one asymptotically linear solution with asymptotic slope $A$ implies that \eqref{tmt} is satisfied for every $M$, with $0 < M < A$. 
\end{remark}

\begin{theorem}
\label{th9}
Let $\sigma$ be right-continuous and nondecreasing on $I$, $F: I\times \mathbb{R}^+\to \mathbb{R^+}$ be continuous and positive in $I\times \mathbb{R}^+$. In addition, let $F(x,y)$ be nondecreasing for every $y>0$, $x\in I$. Then \eqref{nlvs2} has an asymptotically linear solution if and only if \eqref{tmt} holds for some $M>0$.
\end{theorem}

\subsection{Discussion} The previous result extends another result of Nehari [\cite{zn}, Theorem II] to this more general setting. Although we did not exhibit ``Stieltjes analogs" (i.e., for equations of the form \eqref{nlvs1} or \eqref{nlvs2}) of the results in the first few sections for reasons of length, we do not foresee any difficulties in their respective formulations and proofs. In this vein a Stieltjes analog of Theorem~\ref{th6} is readily available, the only major difference being the definition of the space which in this case is $L^{\infty}(I)$. The result is stated next and we leave the proof to the reader.

\begin{theorem}
\label{th10}
Let $f \in L^{\infty}(I)$, $\sigma$ be right-continuous and non-decreasing on $I$, and suppose that the nonlinearity $F:I \times \mathbb{R}\to \mathbb{R}^+$ in 
\begin{equation}
\label{vvv}
y(x) = f(x) - \int_{x}^{\infty} (t-x)\,F(t,y(t))\, d\sigma(t), \quad x\geq x_0
\end{equation}
is continuous on this domain, that $F(x, y)$ is nondecreasing in $y$ for every $x \in I$, $y>0$ and for every $M > 0$, $$\int_{0}^{\infty}t\,F(t, M)\,d\sigma(t) < \infty.$$
In addition, we assume that for every $y, z \in \mathbb{R}$ and every $x \in I$,
$$ | F(x,y)-F(x,z)| \leq k(x)|y-z|$$ where 
$$\int_{x_0}^{\infty}tk(t)\,d\sigma(t) < 1.$$
Then \eqref{vvv} has a solution $y \in L^{\infty}(I)$ if and only if there are two functions $u, v \in L^{\infty}(I)$ such that $u(x) \leq v(x)$, $x \in I$, and for $x\geq x_0$,
\begin{equation}
\label{uf}
u(x) \leq f(x) - \int_{x}^{\infty}(t-x)F(t,v(t))\,d\sigma(t)
\end{equation}
and 
\begin{equation}
\label{vf}
v(x) \geq f(x) - \int_{x}^{\infty}(t-x)F(t,u(t))\,d\sigma(t)
\end{equation}
\end{theorem}

Finally, we give a result that completely parallels Theorem~\ref{th1} above in this wider setting.\\

Let $f\in L^{\infty}[1,\infty)$ with the usual essential supremum norm, $\|\cdot \|$, satisfy \eqref{ft} for some $\delta > 0.$ Define $Y = \{ u\in L^{\infty}[1,\infty): \|u(x)/f(x)\| < \infty\}.$\\

The subset $X = \{ u\in Y: \|u(x)/f(x)\| \leq 2\},$ is a closed subset of $Y$.  Let $F: [1, \infty) \times \mathbb{R} \to \mathbb{R}$ be continuous (and not necessarily positive), and let $\sigma$ be a right-continuous non-decreasing function defined on $[1, \infty)$. In addition, let

\begin{equation} \label{sf00}
\int_{1}^{\infty} s\, | F(s, 0) |\, d\sigma(s) < \infty.
\end{equation}
With $f$ as above let there exist a function $k: [1,\infty) \to \mathbb{R}^+$ satisfying
\begin{equation} \label{sfk0}
\int_{1}^{\infty}\,s\, |f(s)|\, k(s) \, d\sigma(s)  < \infty.
\end{equation}
We assume the usual Lipschitz condition on $F$ as before, that is, for any $u, v \in \mathbb{R}$,  
\begin{equation} \label{lip010}
| F(x, u) - F(x, v) |  \leq k(x) | u-v|, \quad x \geq 1.
\end{equation}

For such functions $F, k, f, \sigma$ satisfying \eqref{ft}, \eqref{sf00}, \eqref{sfk0} and \eqref{lip010} we consider the ``forced" nonlinear equation defined by, for $y \in X$,
\begin{equation}
\label{fnon}
Ty(x) = f(x) - \int_{x}^{\infty} F(t,y(t))\, d\sigma(t), \quad x\geq a.
\end{equation}
on the interval $I = [a, \infty)$ where $a$ is chosen so large that $a \geq 1$ and for $x \geq a$,
\begin{equation} \label{t000}
 {\rm max} \left \{ \int_{x}^{\infty} \,(s-x)\, |f(s)|\, k(s) \, d\sigma(s) ,  \int_{x}^{\infty} \,(s-x)\, |F(s,0)| \, d\sigma(s) \right \} \leq \frac{\delta}{4}.
\end{equation}
 Fix such an $a$ for the next result.
\begin{theorem} 
\label{th15} 
Let $f, F, k, \sigma$ defined above satisfy \eqref{ft}, \eqref{sf00}, \eqref{sfk0}, \eqref{lip010} and \eqref{t000}. Then the operator $T$ has a unique fixed point in $X$, and this point corresponds to a solution of the integral equation 
$$y(x) = f(x) - \int_{x}^{\infty} F(t,y(t))\, d\sigma(t), \quad x\geq a.$$
such that $y \in X$ and $y(x) \sim f(x)$ as $x \to \infty$. 
\end{theorem}

\begin{remark} If $f $ is, in addition, absolutely continuous on $[a,\infty)$, then so is $y$, in which case its right derivative satisfies \eqref{nlvs2} for every $x \geq a$.
\end{remark}

\section{Applications to differential and difference equations}

The main reason for the developments of the previous sections to Volterra-Stieltjes integral and integro-differential equations of the form \eqref{nlvs1}, \eqref{nlvs2} is that this wider framework can be used as a tool for unifying discrete and continuous phenomena such as differential equations and difference equations (or recurrence relations). This approach was emphasized by Atkinson \cite{fvab}, H\"{o}nig \cite{csh}, Mingarelli \cite{abmb} and Mingarelli-Halvorsen \cite{mhb} among the earliest such textual sources. See these texts for basic terminology and other examples of theorems in this wider framework along with their developments to discrete phenomena. Although such generalizations seem to be academic at best, their main thrust lies in their applicability to cases that are not ``continuous" as we see below.\\

The simplest of all applications of the results in Section~\ref{vse} is to differential equations of the second order, linear or not. This is accomplished by choosing $\sigma(t)=t$ throughout that section. The correponding results for ordinary differential equations then arise as corollaries of the results therein. Thus, as pointed out in that section the various theorems therein, some even new for the case of ordinary differential equations, extend essential results in nonlinear theory due to Atkinson, Nehari, Coffman and Wong, etc. to this wider framework. \\

In order to derive results for equations other than ordinary differential equations we can choose $\sigma(t)$ to be a function that is part step-function and part absolutely continuous, or even all step-function or by the same token, all absolutely continuous. The three different choices lead to three intrinsically different kinds of equations. \\

\subsection{The case of three-term recurrence relations}

In order to derive the special results in this case, we appeal to the methods described in [\cite{abmb}, Chapter 1]. Thus, starting from  any infinite sequence of real numbers $\{b_n\}_{n=0}^{\infty}$ we produce an absolutely continuous function $b: \mathbb{N}\to \mathbb{R}$ by simply joining the various points $(n, b_n)$, $n =0,1,2,\ldots$ in the plane by a line segment. The resulting polygonal curve is clearly locally absolutely continuous on its domain (we call this curve the {\it polygonal extension} of the the sequence of points to a curve). Next, we define a right-continuous step-function (or simple function)  by defining its jumps to be at the integers (or any other suitable countable set, [\cite{abmb}, xi]) of magnitude $\sigma(n)-\sigma(n-0) = - b_n$, for $n\geq 0$ (so $\sigma(t) = {\rm constant}$ in between any two consecutive integers). Defining $F(x,y):=y$ for simplicity of exposition, we can show that (see [\cite{abmb}, pp.12-15]) the solution $y(x)$ of the equation \eqref{nlvs2} with right-derivatives has the property that 
\begin{equation*}
{{\Delta}^2}{y_{n-1}}+b_n\,y_n = 0, \quad n\in \mathbb{N},
\end{equation*}
where $y(n)=y_n$ for every $n$, and $\Delta$ is the forward difference operator defined here classically by $\Delta y_{n-1}= y_n - y_{n-1}.$  No more generality is gained by looking at the three-term recurrence relation in standard form, that is,
\begin{equation}
\label{3trr}c_ny_{n+1}+c_{n-1}y_{n-1} + b_ny_n = 0,\quad n\in \mathbb{N},
\end{equation}
where $c_n \neq 0$ for every $n$. The change of dependent variable $y_n = \alpha_nz_n$ where the $\alpha_n$ satisfy the recurrence relation $\alpha_{n+1} = \{c_{n-1}/c_n\}\alpha_{n-1}$, $n\in \mathbb{N}$, brings \eqref{3trr} into the form 
$${{\Delta}^2}{z_{n-1}}+\beta_n\,z_n = 0, \quad n\in \mathbb{N},$$ for some appropriately defined sequence $\beta_n$. Conversely, every such second order linear difference equation  is equivalent to a three term recurrence relation of the form \eqref{3trr} with $y_n=z_n$, $c_n=1$ and $b_n = \beta_n-2$.\\

If $F$ is defined generically as in Section~\ref{vse} then the same choice of the step-function $\sigma$ in \eqref{nlvs2} produces the the second order difference equation
\begin{equation}
\label{2bdef}
{{\Delta}^2}{y_{n-1}}+b_n\,F(n,y_n) = 0, \quad n\in \mathbb{N}.
\end{equation}
The pure nonlinear difference equation 
\begin{equation}
\label{2odef}
{{\Delta}^2}{y_{n-1}}+F(n,y_n) = 0, \quad n\in \mathbb{N}.
\end{equation}
is obtained by setting the $b_n=1$ and defining the resulting step-function $\sigma$ as above. \\

Conversely, starting with any nonlinear difference equation of the form \eqref{2odef} we can produce a Volterra-Stieltjes integro-differential equation of the form \eqref{nlvs2} by ``extending" the domain of this discrete solution $y_n$ to a half axis by joining the points $(n,y_n)$ by line segments. Call this new function $y(x)$. Define the step-function $\sigma$ by jumps of magnitude $\sigma(n)-\sigma(n-0) = - 1$ and right-continuity, and $F(x,y)$, the polygonal extension of the sequence  $F(n, y_n)$ to an absolutely continuous function $F(x,y)$ (obtained by joining the points $(n,y_n,F(n,y_n))$, $(n+1,y_{n+1},F(n+1,y_{n+1}))$, $n\in \mathbb{N}$,  by a line segment). In this case, the Riemann-Stieltjes integral appearing in \eqref{nlvs2} exists for each $x$. The resulting function $y(x)$ is locally absolutely continuous and its right-derivative exists at every point and is locally of bounded variation on the half-axis.  It can be shown that this new function $y(x)$, now satisfies \eqref{nlvs2} with right-derivatives. If more smoothness is required on the function $F$ we can use interpolating polynomials in $\mathbb{R}^3$ {\it in lieu} of the polygonal extension$\ldots$. This duality between equations of the form \eqref{nlvs2} and \eqref{2odef} underlines the importance of this approach.\\

With these facts in hand we formulate the recurrence relation corollary of Theorem~\ref{th5} above.\\

\begin{theorem}
\label{th11}
Let $F:I\times \mathbb{R}\to \mathbb{R}^+$ with values $F(x,y)$,  be continuous on this domain, nondecreasing in its second variable for every $x\in I$ and assume that for some $M>1$ and for some real sequence $\{b_n\}_{n=0}^{\infty}$, we have
$$\sum_{n=0}^{\infty}F(n, Mn)\, |b_n| < \infty.$$
Then \eqref{2bdef} has asymptotically linear solutions, that is solutions of the form $y_n \sim An+B$ as $n \to \infty$ for some constants $A,B$.
\end{theorem}

Another such consequence is a discrete analog of Theorem~\ref{th6}.

\begin{theorem}
\label{th12}
Let $X=\{y\in \mathcal{C}(I): 0\leq y(x)\leq M,\,\,\,\,\, x\in I\}$, where $M>0$ is given and fixed. Let $F: I \times\mathbb{R}^{+}\longrightarrow\mathbb{R}^{+}$ be continuous on this domain, and $\{b_n\}_{n=0}^{\infty}$ a given non-negative sequence such that
\begin{itemize}
\item[$(a)$]$\displaystyle\sum_{n=0}^{\infty}n\,b_n\,F(n,y(n))\leq M,$ \quad{for all $y\in X$,} 
\end{itemize}
\begin{itemize}
\item[$(b)$]$|F(x,u)-F(x,v)|\leq k(x)|u-v|$, \quad{$x \in I$, $u, v \in 
\mathbb{R^+}$}
\end{itemize}
where $k: I \to \mathbb{R}^{+}$ is continuous and for $k(n):=k_n$,
\begin{itemize}
\item[$(c)$]$\displaystyle\sum_{n=0}^{\infty}n\,k_n\,b_n<\infty$.
\end{itemize}
Then the difference equation \eqref{2bdef} has a monotone increasing solution $y_n$ satisfying $0 \leq y_n \leq M$ for each $n$, and $y_n \to M$ as $n \to \infty$.
\end{theorem}

Finally, we formulate a version of Nehari's theorem [\cite{zn}, Theorem I] for second order difference equations as a result of our investigations. We leave the proof to the reader (note that we use $b_n=1$ in this case).

\begin{theorem}
\label{th13}
Let $F: I \times\mathbb{R}^{+}\longrightarrow\mathbb{R}^{+}$ be continuous on this domain and non-decreasing in its second variable (i.e., $F(x, y)$ is nondecreasing in $y$ for $y>0$, for each $x\in I$). Then \eqref{2odef} has bounded eventually positive solutions if and only if 
$$\sum_{n=0}^{\infty}n\,F(n, M) < \infty$$
holds for some $M>0$.
\end{theorem}

This should convince the reader that difference equation analogs of Lemma~\ref{lem1}, Corollary~\ref{cor001}, Corollary~\ref{cor02},Theorem~\ref{th8}, Theorem~\ref{th15} can be formulated without undue difficulty and their proof is simply a consequence of the results in the previous section with the necessary choices of functions as detailed above.\\ 

Next, we note that equations {\it intermediate} between difference and differential equations are also included in our framework of equations of the form \eqref{nlvs2}. That is, we can assume that our function $\sigma$ consists of a discrete part and a part that is possibly continuous and of bounded variation (but not necessarily absolutely continuous). Indeed, on $I=[0,\infty)$ for a given $p>0$ we define $\sigma (t)$ by its jumps on $(0,p]$, so that $\sigma (n) -\sigma(n-1) = -b_n$, for $n=0,1,2,\ldots,p$ where $b_n$ is a given arbitrary sequence and $\sigma$ is right-continuous at its jumps. Let $\sigma (t) := h(t)$ where $h$ is a fixed function, right-continuous and locally of bounded variation on $[p, \infty)$. In the framework of these equations, Nehari's theorem takes the following form:\\

\begin{theorem}
\label{th13}
Let $F: I \times\mathbb{R}^{+}\longrightarrow\mathbb{R}^{+}$ be continuous on this domain and non-decreasing in its second variable (i.e., $F(x, y)$ is nondecreasing in $y$ for $y>0$, for each $x\in I$). Then the integro-differential-difference equation of Stieltjes type,
\begin{equation}
\label{iddi}
y^{\prime}(x) = y^{\prime}(0) - \sum_{n=0}^{p}F(n, y(n))b_n -\int_{p}^{x}F(t,y(t))\,dh(t)
\end{equation}
for $x > p$, has bounded eventually positive solutions if and only if 
$$\int_{p}^{\infty}t\,F(t, M)\,dh(t) < \infty$$
holds for some $M>0$.
\end{theorem}

\subsection{Discussion} A solution $y$ of our equation \eqref{iddi} above is a polygonal curve whenever $0<x<p$ (since the integral term is absent in \eqref{iddi}) while for $x> p$ it is an absolutely continuous curve locally of bounded variation. Thus the values $y(n) := y_n$ actually satisfy a second order difference equation for small $x$ ($x<p$) while for large $x$ ($x>p$) this $y(x)$ is the solution of a pure integral equation of Volterra-Stieltjes type along with some discrete parts (as seen in \eqref{iddi}). The special case $h(t)=t$ is clearly included in this discussion. For this choice, \eqref{iddi} takes the form
\begin{equation}
\label{iddie}
y^{\prime}(x) = y^{\prime}(0) - \sum_{n=0}^{p}F(n, y(n))b_n -\int_{p}^{x}F(t,y(t))\,dt,
\end{equation}
``almost" a second order differential equation except for the interface conditions at a prescribed set of points in $[0,p]$. Under the usual conditions on $F$ as required by Theorem~\ref{th7}, \eqref{iddie} will have eventually positive solutions if and only if 
$$\int_{p}^{\infty}t\,F(t, M)\,dt < \infty$$
holds for some $M>0$ (which is precisely Nehari's necessary and sufficient criterion for second order nonlinear differential equations). For this choice of $\sigma$ this result is to be expected, in some sense, since we are dealing with large $x$ anyhow and so the equation \eqref{iddie} behaves very much like a differential equation. However, we could {\it spread} the discrete part all over the interval $I$ in which case this argument is no longer tenable, as it is {\it a priori} conceivable that oscillations may occur therein (but cannot by Theorem~\ref{th7}).

\section{Proofs}

\begin{proof} (Theorem~\ref{th}) 
We note that $X$ is a closed subset of  the Banach space $Y$ above. This is most readily seen by writing the space $X$ as $X = \{ u \in Y |:  0 \leq \frac{u(t)}{at+b} \leq 1,\text{ for all } t \geq 0\}$ and applying standard arguments. In addition, it is easy to see that $X$ is convex.
Now we define a map $T$ on $X$ by setting
\begin{equation} \label{map}
(Tu)(x) = ax+b - \int_{x}^{\infty} (t-x) \, F(t, u(t))\, dt  
\end{equation}
for $u \in X$. Note that the right-side of \eqref{map} converges 
for each $x \geq 0$, because of \eqref{atk}. Indeed, for $u \in X$, $x \geq 0$, 
\begin{equation} \label{bound}
0 \leq  \int_{x}^{\infty}(t-x)F(t, u(t))\, dt \leq \int_{0}^{\infty} t\, 
F(t, u(t))\, dt \leq b,
\end{equation}
as $F(t, u(t)) \geq 0$ for such $u$ (which implies that $(Tu)(x) \leq ax+b$) 
and the indefinite integral is a non-increasing function of $x$ on 
$[0, \infty)$. Since $a \geq 0$, we get that $(Tu)(x) \geq 0$ for any $x \geq 0$. On the other hand, it is easy to see that for $u \in X$, $Tu$ is a continuous function on $[0, \infty)$. So, $TX \subseteq X$.

Next, we prove that $T$ is a continuous map on $X$. For $u, v \in X$, 
\begin{align*}
|(Tu)(x) - (Tv)(x) | 
& \leq   \int_{x}^{\infty} (t-x) | F(t, u(t)) - F(t, v(t)) | \, dt  \\
& \leq   \int_{x}^{\infty} (t-x)  k(t) | u(t) - v(t) |\, dt  \\
& \leq  \| u - v \|_{Y} \, \int_{0}^{\infty} t\, k(t)\,(at+b)\, dt,
\end{align*}
where we have used \eqref{lip} and the fact that 
$\int_{x}^{\infty}(t-x)k(t)(at+b)\, dt$ is a non-increasing function of $x$ for 
$x \in [0, \infty)$, since $k(t)(at+b) \geq 0$. It follows that for $x \geq 0$, 
\begin{eqnarray}\label{cont}
| \Psi (Tu)(x) - \Psi (Tv)(x) | & \leq & \frac{1}{b}\, \| u - v \|_{Y} \, \max\{a,b\}\, \int_{0}^{\infty} t\,(t+1)\, k(t)\, dt,
\end{eqnarray}
from which we conclude that
$$
\| Tu - Tv \|_{Y} \leq \alpha \, \| u - v \|_{Y},
$$ 
where $\alpha < \infty$ on account of \eqref{teek} and \eqref{tktt}. It follows that $T$ is continuous on $X$.

Next, we show that $TX$ is compact, that is, $T$ sends bounded subsets of $X$ onto relatively compact subsets. For $M$ a subset 
of $X$ we have that to prove that $TM$ is relatively compact. By virtue of the isometry $\Psi$, this is equivalent to proving that 
$\Psi (T(M))$ is relatively compact. To this end, we use the measure of noncompactness on $BC(\mathbb{R}^+)$ defined for $A \in BC(\mathbb{R}^+)$ by
\[
\mu (A)= \lim_{L\to \infty}\left( \lim_{\varepsilon \to 0} w^L (A,\varepsilon)\right) + \limsup_{t\to\infty}\, {\rm diam}\, A(t),
\]
see [\cite{bg}, Theorem 9.1.1(d), p.46], where
$${\rm diam}\, A(t) = {\rm sup} \{|x(t)-y(t)| : x, y \in A\},$$
and 
$$w^L (A,\varepsilon) = {\rm sup} \{ w^L(x,\varepsilon) : x \in A\},$$ with $$w^L(x,\varepsilon) = {\rm sup} \{|x(t)-x(s)| : t, s \in [0,L], |t-s| \leq \varepsilon\}.$$

We fix $\varepsilon >0$, $L>0$, $u\in M\subset X$ and $t_1,t_2\in\mathbb{R}^+$ with $t_2-t_1\leq \varepsilon$ and, 
without loss of generality $t_2>t_1$. Then 
\[
\begin{array}{l}
|\Psi(Tu)(t_2)-\Psi(Tu)(t_1)|= \left|\displaystyle\int_{t_2}^{\infty}\frac{ (s-t_2)F(s,u(s))}{at_2 +b}ds- 
\displaystyle\ \int_{t_1}^{\infty}\frac{ (s-t_1)F(s,u(s))}{at_1 +b}ds \right| \medskip \\
\text{\ \ }\leq \left|\displaystyle\ \int_{t_2}^{\infty}\left[\frac{(s-t_2)}{at_2+b}-\displaystyle\ 
\frac{(s-t_1)}{at_1 +b}\right]F(s,u(s))ds -
\displaystyle\ \int_{t_1}^{t_2} \frac{(s-t_1)}{at_1 +b}F(s,u(s))ds\right| \medskip\\
\text{\ \ } \leq \displaystyle\ \int_{t_2}^{\infty}\frac{(as +b)(t_2-t_1)}{(at_2+b)(at_1+b)}F(s,u(s))ds + \displaystyle\ 
\int_{t_1}^{t_2}\frac{(s-t_1)}{at_1 +b}F(s,u(s))ds 
\medskip\\
\text{\ \ }\leq  \displaystyle\ \frac{\varepsilon}{b^2}\displaystyle\ \int_{t_2}^{\infty}(as +b) F(s,u(s))ds + \displaystyle\ 
\frac{\varepsilon}{b}\int_{t_1}^{t_2}F(s,u(s))ds  \medskip\\
\text{\ \ }\leq  \displaystyle\ \frac{a \varepsilon}{b} + \displaystyle\
\frac{\varepsilon}{b}\int_{t_1}^{\infty}F(s,u(s))ds,
\end{array}
\]
by \eqref{atk}, since $F \geq 0$ for $u \in M$. Combining these estimates we deduce that
\begin{eqnarray}
|\Psi(Tu)(t_2)-\Psi(Tu)(t_1)| & \leq &  \displaystyle\ \frac{a \varepsilon}{b} + \displaystyle\
\frac{\varepsilon}{b}\int_{t_1}^{\infty}F(s,u(s))ds, \nonumber \\
& \leq &  \displaystyle\ \frac{a \varepsilon}{b} + \frac{\varepsilon}{b}\left \{ \int_{t_1}^{t_2}F(s,u(s))ds + \int_{t_2}^{\infty}F(s,u(s))ds \right \},\nonumber \\ 
& \equiv & \label{28} \displaystyle\ \frac{a \varepsilon}{b} + \frac{\varepsilon}{b}\left \{ I_1 +I_2 \right \}.
\end{eqnarray}
We estimate the two integral quantities in \eqref{28} in turn. This said, use of  \eqref{lip} and \eqref{atk}  for $u \in M$ gives
\begin{eqnarray}
I_2&\leq & \int_{t_2}^{\infty} |F(s,u(s))-F(s,0)|ds+ \int_{t_2}^{\infty} F(s,0)ds \nonumber \\
&\leq &  \int_{0}^{\infty} k(s)u(s)ds + \int_{0}^{1}F(s,0)ds +  \int_{1}^{\infty}sF(s,0)ds \label{28a} \\
&\leq & \|u\|_Y \int_0^{\infty}k(s)(as +b)ds  + \displaystyle \sup_{s \in [0,1]}F(s,0) +\int_{0}^{\infty}sF(s,0)ds \nonumber  \\
&\leq & C_2 \equiv \int_0^{\infty} k(s)(as +b)ds +  \displaystyle \sup_{s \in [0,1]}F(s,0) + b \label{28b}
\end{eqnarray}
(since, for $u\in M$, $\|u\|_Y \leq 1$), where $C_2$ is finite and independent of $\varepsilon$, because of \eqref{tktt} and \eqref{kt}. 
On the other hand, arguing as in \eqref{28a}-\eqref{28b}, we obtain
\begin{eqnarray}
I_1&\leq  & \int_{t_1}^{t_2} |F(s,u(s))-F(s,0)|ds+ \int_{t_1}^{t_2} F(s,0)ds \nonumber \\
&\leq &  \int_{0}^{\infty} k(s)u(s)ds + \int_{0}^{\infty}F(s,0)ds \nonumber \\
&\leq & C_2 \label{29a}.
\end{eqnarray}
Combining \eqref{29a}, \eqref{28b} and \eqref{28} we get finally,
\begin{eqnarray}
|\Psi(Tu)(t_2)-\Psi(Tu)(t_1)| & \leq & \displaystyle\ \frac{a \varepsilon}{b} + \frac{2C_2 \varepsilon}{b}. \label{30}
\end{eqnarray}
Passing to the supremum over $u \in M \subset X$ we find that
\[
\lim_{\varepsilon \to 0} w^{L}(\Psi (TM),\varepsilon)=0
\]
Since $L>0$ is arbitrary we deduce that 
\[
\lim_{L\to \infty} \left(\lim_{\varepsilon \to 0}w^L (\Psi (TM),\varepsilon)\right)=0.
\]

In order to complete the proof we need to analyse the term related to the diameter. Taking $u, v\in M$ and $t\in  \mathbb{R}^+$ then, proceeding as in the continuity argument above leading to \eqref{cont}, we see that 
\begin{eqnarray*}
|\Psi (Tu)(t)-\Psi (Tv)(t)| &\leq & \displaystyle \frac{\|u-v\|_{Y}}{b}\displaystyle \int_t^{\infty} s(as+b)k(s)ds \\
& \leq & \displaystyle \frac{\|u-v\|_Y}{b}\left[a \displaystyle \int_t^{\infty} s^2 k(s)ds +b\displaystyle \int_t^{\infty}sk(s)ds\right] \\
& \leq & \frac{2}{b} \left[a \displaystyle \int_t^{\infty} s^2 k(s)ds +b\displaystyle \int_t^{\infty} sk(s)ds\right],
\end{eqnarray*}
since $\|u-v\|_Y \leq 2$ for $u, v \in M$. Consequently 
\[ 
\limsup_{t\to \infty}\, {\rm diam}\, \Psi (TM)(t)\, \leq \frac{2}{b}\limsup_{t\to \infty}\left[ a\int_t^{\infty} s^2 k(s)ds +b\int_t^{\infty} sk(s)ds \right]=0,
\]
on account of \eqref{teek} and Remark~\ref{uno}. Thus, $\mu (\Psi (TM))=0$. This fact tells us that $\Psi (TM)$ is relatively compact in $BC(\mathbb{R}^+)$ and, since $\Psi$ is an isometry, $TM$ is relatively compact in $Y$. Hence, $T$ is compact, and so Schauder's theorem gives the existence of a fixed point $u \in X$ for $T$. This fixed point is necessarily a solution of \eqref{non0} asymptotic to the line $ax+b$ as $x \to \infty$. This completes the proof. 
\end{proof}

\begin{proof} ({Theorem~\ref{th1}})  Let $Y$ be the Banach space defined in \eqref{spaceY} with the norm \eqref{normY}, where we replace the interval $[1,\infty)$ by $x_0,\infty)$. Let $X$ be the closed subset defined by
$X =\{u\in Y: {\sup_{x \geq x_0}\{|u(x)|/|f(x)|\}\leq 2} \}.$ Define a map $T$ on $X$ by $u \in X$, 
\begin{equation}
\label{tux}
Tu(x) = f(x) - \int_{x}^{\infty}{(s-x)F(s,u(s))\,ds},\ \  x\geq x_0.
\end{equation} Clearly, for $u \in C(I)$ we have, because of our assumptions on $F$, $Tu \in C(I)$. In addition, \eqref{lip01} gives that for $u \in X$, $|F(s,u(s))| \leq k(s)|u(s)| + |F(s,0)|$. Combining this with \eqref{tux}, dividing \eqref{tux} throughout by $f(x)$ a simple estimation gives that for $ x \geq x_0$, 
\begin{eqnarray*}
\bigg | \frac{Tu(x)}{f(x)}\bigg | \leq 1 + \frac{1}{|f(x)|}\int_{x}^{\infty}(s-x)|k(s)f(s)|\bigg |\frac{u(s)}{f(s)}\bigg | \, ds \,+\\ \frac{1}{|f(x)|}\int_{x}^{\infty}(s-x)|F(s,0)|\, ds. \\ 
\end{eqnarray*} Now, use of \eqref{ft} shows that 
\begin{eqnarray}
\label{Tuf}
\bigg | \frac{Tu(x)}{f(x)}\bigg | \leq 1 + \bigg \| \frac{u}{f}\bigg \| \frac{1}{\delta}\int_{x}^{\infty}(s-x)|k(s)f(s)|\, ds \,+ \nonumber \\ \frac{1}{\delta}\int_{x}^{\infty}(s-x)|F(s,0)|\, ds,
\end{eqnarray} which, since $u \in X$ and \eqref{t0} is enforced, furnishes the bound 
\begin{eqnarray*}
\bigg \| \frac{Tu(x)}{f(x)}\bigg \| \leq 1 + 2 \frac{1}{\delta}\frac{\delta}{4}\,+ \frac{1}{\delta}\frac{\delta}{4} = \frac{7}{4} < 2.\end{eqnarray*}  Thus $T$ is a self-map on $X$.  In order to show that $T$ is a contraction on $X$, consider the simple estimate derived from \eqref{tux}, namely, for $x \geq x_0$,
\begin{eqnarray*}
\bigg | \frac{Tu(x) - Tv(x)}{f(x)} \bigg | &\leq& \frac{1}{|f(x)|} \int_{x}^{\infty}{(s-x)|F(s,u(s))-F(s,v(s))|\, ds},\\
&\leq& \frac{1}{|f(x)|} \int_{x}^{\infty}{(s-x)k(s)|u(s)-v(s)|\, ds},\quad ({\rm by}\, \eqref{lip01} ) \\
&\leq& \frac{1}{\delta}\|u-v\| \int_{x}^{\infty}{(s-x)k(s)|f(s)| \, ds}, \quad ({\rm by}\, \eqref{ft} ) .
\end{eqnarray*}
Since the last display is valid for every $x \geq x_0$ it follows from \eqref{t0} that,
$$\| Tu - Tv\| \leq (1/4)\, \|u-v\|,$$
so that $T$ is a contraction on $X$. It is easily seen that the subsequent fixed point, say $u(x)$, obtained by applying the classical fixed point theorem of Banach, is a solution of \eqref{non0} satisfying the conclusion (2) stated in the theorem, since $u \in X$. On the other hand, since our fixed point $u$ satisfies \eqref{tux}, we have
$$\frac{u(x)}{f(x)} = 1 - \frac{1}{f(x)}\int_{x}^{\infty} (s-x)F(s,u(s))\, ds.$$ An estimation of this integral similar to the one leading to the right-side of \eqref{Tuf} gives that $$\lim_{x \to \infty} \frac{1}{f(x)}\int_{x}^{\infty} (s-x)F(s,u(s))\, ds = 0,$$ on account of the finiteness of all the integrals involved. This shows that $u(x) \sim f(x)$ as $x \to \infty$.
\end{proof}

\begin{proof} (Theorem~\ref{th2}) Note that $X$ is a closed subset of the Banach space $BC(\mathbb{R}^+)$. For $u \in X$ we define a map $T$ by setting $$Tu(x) = f(x) - \int_{x}^{\infty}(t-x)F(t, u(t))\, dt, \quad x \geq 0.$$ Then for $u \in X$ it is clear that $Tu \in C(\mathbb{R}^+)$ and since $F(t,u(t)) \geq 0$ for such $u$ and all $t \geq 0$, we have $|Tu(x)| \leq \|f\|_{\infty} + b$, for every $x \geq 0$, where we have used the fact the integral in question is a non-increasing function of $x$ for all $x \geq 0$. Hence $T$ is a self-map on the ball $X$. Finally, an argument similar to the corresponding one in Theorem~\ref{th} gives that $T$ is a contraction on $X$ provided there holds \eqref{tkt}. This completes the proof.
\end{proof}

\begin{proof} (Theorem~\ref{th4}): The necessity is simple. If $y$ is such a solution then set $u=v=y$ throughout. For the sufficiency we appeal, as usual, to a fixed point theorem. Consider the space $BC(I)$ of (uniformly) bounded continuous functions on $I$ with the uniform norm. Since $u, v$ are uniformly bounded by hypothesis, the subset $X$ defined by $$X=\{ y \in BC(I) : u(x) \leq y(x) \leq v(x), x \in I\}$$ with the induced metric, is complete. Define a map $T$ on $X$ by the usual
\begin{equation*}
Ty(x) = f(x) - \int_{x}^{\infty}(t-x)F(t,y(t))\,dt, \quad x \in I
\end{equation*}
for $y \in X$. Since $y \in X$, then $y \in L^{\infty}(I)$; it follows from hypotheses (2) and (3) that the integral on the right is finite for every $x \in I$ and this defines a continuous function that is uniformly bounded on $I$. Thus, $Ty$ is continuous and uniformly bounded on $I$, since $f$ is. Thus, $T$ is well-defined. On the other hand, by hypothesis (2), $F(t,u(t)) \leq F(t,y(t)) \leq F(t,v(t))$ for $t \in I$; it follows that, for $x \in I$,
\begin{eqnarray*}
Ty(x) &\leq& f(x) - \int_{x}^{\infty}(t-x)F(t,u(t))\,dt \leq v(x)\\
& \geq & f(x) - \int_{x}^{\infty}(t-x)F(t,v(t))\,dt \geq u(x)
\end{eqnarray*}
 where we have used assumptions (b) and (c) in order to estimate the integrals. Thus $T$ is self-map on $X$. That $T$ is a contraction on $X$ follows the usual route. Briefly, for $x \in I$, $y, z \in X$,
\begin{eqnarray*}
|Ty(x) - Tz(x)| &\leq& \int_{x}^{\infty}(t-x)|F(t,y(t))-F(t,z(t))|\, dt\\
&\leq & \int_{x}^{\infty}(t-x)k(t)|y(t)-z(t)|\, dt\\
& \leq & \left (\int_{x_0}^{\infty}tk(t)\,dt\right ) \,\, \|y-z\|_{\infty}
\end{eqnarray*}
and so,
$$\|Ty - Tz\|_{\infty} \leq  \alpha \|y-z\|_{\infty},$$
where $\alpha < 1$ is the integral in question (cf., assumption (5)).
\end{proof}

\begin{proof} (Theorem~\ref{cor01}) This is clear since we can integrate the differential equation \eqref{non0} twice to obtain \eqref{ie1} and conversely, if we know that its solution is $L^{\infty}$, we can differentiate \eqref{ie1} twice to recover \eqref{non0}. The result follows from an application of the theorem.
\end{proof}

\begin{proof}(Theorem~\ref{thx}) We integrate the inequalities twice over the half line to obtain both \eqref{uf} and \eqref{vf}. An application of Theorem~\ref{th4} gives that \eqref{ie1} has a solution $y(x) \sim f(x)$, as $x \to \infty$. But the right side of \eqref{ie1} is twice differentiable, consequently so is $y(x)$, that is \eqref{pert} is satisfied.
\end{proof}

\begin{proof}(Theorem~\ref{atm}) First, we show that solutions of \eqref{non2} exist on the half-line, $I$. Introduce the usual energy functional $E(x)$ on solutions of \eqref{non2} by 
\begin{equation}\label{ener} E(x) = \frac{1}{2}{y^{\prime}}^2 +  \int_{0}^{y}\eta G(x, \eta)\, d\eta \equiv \frac{1}{2}{y^{\prime}}^2 + \mathcal{I}(x,y),
\end{equation}
where $\mathcal{I}_{x}(x,y) \leq 0$ by hypothesis. A glance at \eqref{non2} shows that $$E^{\prime}  = y^{\prime}g + \mathcal{I}_{x} \leq y^{\prime}g \leq g \sqrt{2E}.$$ So, $E^{\prime}E^{-1/2} \leq \sqrt{2}g$ whenever $E>0$. It follows that if the solution $y(x)$ exists for $x \in [a,b]$, $a \geq x_0$ and, at the same time, $E(x) > 0$ for such $x$,  then
\begin{equation}
\label{ener2} \sqrt{E(b)} \leq \sqrt{E(a)} + \frac{\sqrt{2}}{2}\int_{a}^{b}g(t)\,dt.
\end{equation}
Of course, \eqref{ener2} is also true for any interval $[a,b]$ in which the solution exists. For such an interval we have from \eqref{ener2} 
\begin{equation}
\label{ener3} |y^{\prime}(b)| \leq \sqrt{2E(a)} + \int_{a}^{b}g(t)\,dt,
\end{equation}
so, if the solution exists on an interval $[a,b)$ then it can be continued to $x=b$ and thus to a right-neighborhood of $b$. Thus, we see that for any $x \geq x_0$ a solution can be continued throughout $I$.

We now claim that for a given solution $y$ of \eqref{non2} there is an $X$ (depending on $y$) such that we cannot have for 
\begin{equation}
\label{yayb}b > a \geq X, y(a)=y(b)=0, y(x) > 0, x \in (a,b).
\end{equation}
Note that by \eqref{tgt} and \eqref{ener3} we can suppose that $X$ is such that $X > 0$ and such that for some $K>0$ we have 
\begin{equation}
\label{yKx} |y(x)| < Kx, \quad x \geq X.
\end{equation}
This already implies that all solutions are ``sublinear" or cannot grow faster than a linear function. We fix this $K$ and consider the differential equation $$z^{\prime\prime} + G(x, Kx)z =0, \quad x \geq X.$$ Since this is a linear equation it is well known that (e.g., \cite{rb}),  assumption \eqref{tgKt} implies that this equation has a solution $z(x) \to 1$ as $x \to \infty$. We choose $X_0> X$ so that $z(x) > 0$ for all $x \geq X_0$ and let $b>a \geq X_0$. Now, writing $y = wz$ and making use of the equation for $z$, we obtain the second order linear differential equation 
\begin{equation}
\label{wz}
w^{\prime\prime}z+2w^{\prime}z^{\prime} = g+wz\{G(x,Kx)-G(x,y)\}.
\end{equation}
However, \eqref{yayb}, \eqref{yKx} and $G$ non-decreasing in its second variable,  would imply that the right of \eqref{wz} is positive in $(a,b)$ so that $(w^{\prime}z^2)^{\prime} \geq 0$ on $(a,b)$. Indeed, \eqref{yayb} would also force $w(a)=w(b)=0$ and $w(x)>0$ in $(a,b)$. On the other hand, this leads to $w^{\prime}(a) \geq 0$, $w^{\prime}(b) \leq 0$. Since $w^{\prime}z^2$ is non-decreasing, this implies that $w^{\prime}=0$, i.e., $w=0$ in $(a,b)$ resulting in a contradiction. Thus, $y(x)$ is either ultimately positive or it is ultimately non-positive.
\vskip0.15in

Now consider the case where $y(x)$ is ultimately positive. We may suppose (see \eqref{yKx}) that 
\begin{equation}
\label{0yKt}
0 < y(x) < Kx, \quad x \geq X.
\end{equation}
We modify the argument following \eqref{yKx} as follows: Consider the differential equation $$z_1^{\prime\prime}+G(x,y(x))z_1 =0 , \quad x\geq X.$$ As before, the integrability condition \eqref{tgKt} gives that this will have a solution $z_1(x)$ with $z_1(x) \to 1$ as $x \to \infty$. We can define $w_1$ as before by $y=w_1z_1$ and find, as before, that $(w_1^{\prime}{z_1}^2)^{\prime} =gz_1.$ But \eqref{tgt} along with the fact that $w_1^{\prime}{z_1}^2$ is non-decreasing implies that ${w_1}^{\prime}$ tends to a non-negative finite limit at infinity. The possibility that ${w_1}^{\prime}(\infty)>K$ is excluded on account of \eqref{0yKt}. Hence $ 0 \leq {w_1}^{\prime}(\infty) \leq K$. If $A\equiv {w_1}^{\prime}(\infty)>0$  then necessarily $y(x) \sim Ax$ as $x \to \infty.$

The other possibility is that \eqref{0yKt} holds but that ${w_1}^{\prime}(\infty)=0$. In this case, the differential equation for $w_1$ yields
$${w_1}^{\prime}(x) = - \{z_1(x)\}^{-2}\,\int_{x}^{\infty}g(t)z_1(t)\,dt,$$
and since $z_1(x) \to 1$ as $x \to \infty$ we see that 
\begin{equation}
\label{w1}
{w_1}^{\prime}(x) \sim - \int_{x}^{\infty}g(t)\,dt
\end{equation}
as $x \to \infty$. Note that if \eqref{tgt} were false (for $i=1$)  it would follow that (since $g(x) \geq 0$), $w_1(x) \to -\infty$ as $x \to \infty$ and this contradicts the positivity of $y(x)$ for all large $x$. Thus, \eqref{tgt} is actually a necessary condition. On the other hand, the hypothesis \eqref{tgt} implies that $B\equiv w_1(\infty)$ is finite and necessarily non-negative, because of the positivity of $y$, i.e., $y(x) \sim B$ as $x \to \infty$, where $B \geq 0$.
\vskip0.15in
If $y(x)$ is ultimately non-positive, we can take it that $y(x) \leq 0$  for $x \geq X$, and that $y(x) < 0$ on some unbounded subset of $x \geq X$. Since $yG(x,y) \leq 0$ for $y \leq 0$, we get from \eqref{non2} that $y^{\prime\prime} \geq 0$. Applying Lemma 0 in \cite{fva} with $z(x)\equiv -y(x)$ we see that $y^{\prime}(x) \leq 0$ for $x \geq X$. This, in conjunction with the fact that $y(x) < 0$ on some unbounded subset implies that $y(x) < 0$ for all sufficiently large $x$. So, $ y(x) < 0$, $y^{\prime}(x) \leq 0$, $y^{\prime\prime} \geq 0$ for all large $x$ which leads to a counterpart of the positive solutions result.
\end{proof}

\begin{proof} (Theorem~\ref{atm2}) As before we write $F(x,y)=yG(x,y)$ and we proceed as in the proof of Theorem~\ref{atm} up to \eqref{ener3}. The same argument therein gives that solutions all exist on some half-axis. Indeed, since $f^{\prime}(\infty) = \infty$, by assumption, we can use \eqref{ener3} to derive that $|y^{\prime}(x)| \leq (1+\varepsilon)f^{\prime}(x)$ for all sufficiently large $x$. In addition, another integration gives us a similar bound for $y$ in the form $|y(x)| \leq (1+\varepsilon)f(x)$ for all sufficiently large $x$, say $x \geq X_0$. Next, an integration of \eqref{non2} over $[X_0, x)$ and use of \eqref{tgt2} and \eqref{fr} shows that $$y^{\prime}(x) \sim \int_{x_0}^{x}g(t)\,dt, \quad x\to \infty.$$ Finally, one last integration of the preceding equation gives the desired asymptotic estimate.
\end{proof}

\begin{proof}(Theorem~\ref{th5})
Consider the solution $y$ whose initial conditions are $y(a)=a$, $y^{\prime}(a)=M$ with right-derivatives, where $a$ is to be chosen later. Since $y^{\prime}$ is right-continuous on $I$ there is a $b > a$ such that 
\begin{equation}
\label{m1}
\frac{M}{2} < |y^{\prime}(x)| < 2M, \quad{x \in [a,b)}
\end{equation}
Since $y$ is absolutely continuous on $[a,b)$ it follows that
\begin{equation*}
|y(x)| \leq |y(a)| + |\int_{a}^{x}y^{\prime}(t)\, dt | \leq a+ M(x-a)
\end{equation*} 
that is,
\begin{equation}
\label{m2}
|y(x)| < Mx, \quad{x \in [a,b)}
\end{equation}
since $M > 1$. Since $y$ is continuous, \eqref{m2} also holds at $x=b$. It follows from \eqref{nlvs1}, \eqref{nlvs2} that
\begin{equation}
\label{m3}
y^{\prime}(b) = y^{\prime}(a) - \int_{a}^{b}F(t,y(t))\,d\sigma(t)
\end{equation}
i.e.,
\begin{equation}
\label{m4}
| y^{\prime}(b) - M|\leq  \int_{a}^{b}F(t,y(t))\,|d\sigma(t)|.
\end{equation}
On the other hand, $F$ is nondecreasing in its second variable by hypothesis, so \eqref{m2} and assumption (3) together yield
\begin{equation}
\label{m5}
| y^{\prime}(b) - M|\leq  \int_{a}^{\infty}F(t,Mt)\,|d\sigma(t)|.
\end{equation}
Now, by assumption (3) again we can choose (and fix) $a$ so large that $$\int_{a}^{\infty}F(t,Mt)\,|d\sigma(t)| < M/4.$$ Then \eqref{m1} holds for all $b > a$ and, as a result, \eqref{m2} holds for all $x > a$. From this we see that for given $\varepsilon > 0$, we can choose $X$ so large that for any $c> b > X$ we have 
$$\int_{b}^{c}F(t,Mt)\,|d\sigma(t)| < \varepsilon,$$ and a double application of \eqref{m3} and the usual estimates, shows that $$|y^{\prime}(c)-y^{\prime}(b)| <\int_{b}^{c}F(t,Mt)\,|d\sigma(t)| < \varepsilon,$$ for $c > b > X.$ Hence $y^{\prime}(x)$ tends to a limit $L$, say, as $x \to \infty$, and $L \neq 0$ on account of \eqref{m1}. From this it follows that $y(x)/x = L + o(1)$, i.e., $y$ is asymptotically linear as $x \to \infty$.
\end{proof}

\begin{proof}(Theorem~\ref{th6}) It suffices to show that the
 integral equation
\begin{equation}\label{eq1}
y(x)=M-\displaystyle\int_{x}^{\infty}(t-x)\,F(t,y(t))\, d\sigma(t),\end{equation}
has a fixed point in $X$ (since the resulting solution will be absolutely continuous, with a derivative that is locally of bounded variation and satisfying \eqref{nlvs2}).  So, we define the operator $A$ on $X$ by 
$$(Ay)(x)=M-\displaystyle\int_{x}^{\infty}(t-x)F(t,y(t))\,d\sigma(t).$$ Since $\sigma$ is non-decreasing and $F \geq 0$, for $y\in X$, the function defined by $$\int_{x}^{\infty}(t-x)\,F(t,y(t))\,d\sigma(t)$$ is nonincreasing, hence
$$ 0\leq\int_{x}^{\infty}(t-x)\,F(t,y(t))\,d\sigma(t)\leq
\int_{0}^{\infty}t\,F(t,y(t))\,d\sigma(t)\leq M,
$$
by hypothesis (a). Consequently, for $y\in X$ we have y
\begin{equation}\label{0ay}
0\leq (Ay)(x)\leq M,\makebox[1cm]{for}x\in I.
\end{equation}
In order to show that for $y\in X$ then
$Ay\in\mathcal{C}(I)$, we note by Fubini's theorem that since 
$$\int_{x}^{\infty}(t-x)F(t,y(t))\, d\sigma(t) = \int_{x}^{\infty}\int_{t}^{\infty}F(s,y(s))\, d\sigma(s)\, dt,$$
for every $x \in I$ and the integral of a function that is locally of bounded variation is locally absolutely continuous, it is in particular continuous and so, for $y\in X$, we have $Ay\in\mathcal{C}(I)$. This, in combination with \eqref{0ay} shows that $Ay \in X$. Hence, the operator $A$ applies $X$ into itself. \\ 

\noindent Now, we prove that $A$ is continuous on $X$. Indeed,
\begin{eqnarray*}
|(Ay_{n})(x)-(Ay)(x)|&= & \left|\displaystyle\int_{x}^{\infty}(t-x)[F(t,y_{n}(t))-F(t,y(t))]\, d\sigma(t)\right| \\
&\leq &\displaystyle\int_{x}^{\infty}(t-x)|F(t,y_{n}(t))-F(t,y(t))|\, d\sigma(t) \\
&\leq &\|y_{n}-y\|_{\infty}\displaystyle\int_{x}^{\infty}(t-x)k(t)\, d\sigma(t) \\
&\leq &\|y_{n}-y\|_{\infty}\displaystyle\int_{0}^{\infty}t\,k(t)\,d\sigma(t).
\end{eqnarray*}
\noindent It follows that $A$ is continuous on $X$ on account of assumption $(c)$. \\

\noindent The proof that $A$ is compact uses ideas from the theory of measures on non-compactness. First, we introduce some terminology. Let us fix a nonempty bounded subset $X$ of $C[0,a]$. For 
$\varepsilon >0$ and $y \in X$ denote by 
$w(y,\varepsilon)$ the {\it modulus of continuity} of $y$ defined by 
\[
w(y,\varepsilon)= sup \ \{\mid y(t)-y(s)\mid \ :\ t,s \in [0,a],\ \mid t-s\mid \leq \varepsilon\}
\]
Further, let us put 
\begin{eqnarray*}
w(X,\varepsilon)&=&sup\ \{w(y,\varepsilon)\ :\ y\in X\}\\
w_0 (X)&=& \lim_{\varepsilon \rightarrow 0} w(X,\varepsilon),
\end{eqnarray*}
It can be shown (see \cite{11}) that the function $\mu(X)=w_0 (X)$ is a regular measure of noncompactness in the space $C[0,a]$. Now, let $x_{1},x_{2}\in [0,\infty)$ be such that $x_{2}-x_{1}\leq \varepsilon$ and without loss of generality, $x_{1}<x_{2}$. Then 
\begin{eqnarray}
|Ay(x_{2})-Ay(x_{1})|
&=&\left|\int_{x_{2}}^{\infty}(t-x_{2})F(t,y(t))d\sigma(t)-
\int_{x_{2}}^{\infty}(t-x_{1})F(t,y(t))d\sigma(t)+ \right. \nonumber\\
& & + \left. \int_{x_{2}}^{\infty}(t-x_{1})F(t,y(t))d\sigma(t)- 
\int_{x_{1}}^{\infty}(t-x_{1})F(t,y(t))d\sigma(t)\right|\nonumber \\
&\leq &\left|\int_{x_{2}}^{\infty}(x_{1}-x_{2})F(t,y(t))d\sigma(t)+\int_{x_{1}}^{x_{2}}(t-x_{1})F(t,y(t))d\sigma(t)\right|\nonumber\\
&\leq &\int_{x_{1}}^{x_{2}}(t-x_{1})F(t,y(t))d\sigma(t)+
\int_{x_{2}}^{\infty}(x_{2}-x_{1})F(t,y(t))d\sigma(t)\nonumber\\
&\leq & (x_{2}-x_{1})\int_{x_{1}}^{\infty}F(t,y(t))d\sigma(t).\label{mnc0}
\end{eqnarray}
Next, we note that for any $x_0 \geq 0$, $y \in X$,
\begin{eqnarray}
\int_{x_{0}}^{\infty}F(t,y(t))\,d\sigma(t) &\leq &
\int_{0}^{\infty}F(t,y(t))\, d\sigma(t)\nonumber \\
&\leq &
\int_{0}^{1}F(t,y(t))\,d\sigma(t)+\int_{1}^{\infty}t\,F(t,y(t))\,d\sigma(t)\nonumber\\
&\leq &\|F\|_{[0,1]\times [0,M]}(\sigma(1)-\sigma(0))+M,\label{mnc}
\end{eqnarray}
by hypothesis (a) and since $\sigma$ is nondecreasing. Thus, use of \eqref{mnc} and \eqref{mnc0} give us that
$$
w(Ay,\varepsilon)\leq\varepsilon[\|F\|_{[0,1]\times
[0,M]}(\sigma(1)-\sigma(0))+M],
$$
\noindent consequently,
$$
w(AX,\varepsilon)\leq \varepsilon[\|F\|_{[0,1]\times
[0,M]}(\sigma(1)-\sigma(0))+M],
$$
\noindent so that
\begin{equation}
\label{w0a}
w_{0}(AX)=0.
\end{equation}
\noindent Finally, let $y,z\in X$, $x\geq 0$. Then the previous continuity argument also yields the estimate
\begin{equation}
\label{ayxaz}
|(Ay)(x)-(Az)(x)|\leq \|y-z\|_{\infty}\,\int_{x}^{\infty}(t-x)\,k(t)\,dt
\end{equation}
On the other hand, for $y, z\in X$,  $\|y-z\|_{\infty}\leq 2M$ and so
\begin{eqnarray*}
|(Ay)(x)-(Az)(x)|&\leq & 2M\int_{x}^{\infty}(t-x)k(t)\,d\sigma(t)\\
&\leq & 2M\displaystyle\int_{x}^{\infty}t\,k(t)\, d\sigma(t),
\end{eqnarray*}
since $x \in I$. It follows that the diameter of the set $AX$ can be estimated by
$$
diam AX(x)\leq 2M\displaystyle\int_{x}^{\infty}sk(s)d\sigma(s),
$$
and taking the limit as $x\to \infty$, we get
\begin{equation}\label{axx}
\lim_{x\to\infty}diam (AX)(x)=0.
\end{equation}
\noindent Therefore, \eqref{w0a} and \eqref{axx} give us that $AX$ is compact. An application of Schauder's fixed point theorem now gives the desired conclusion.
\end{proof}

\begin{proof} (Corollary~\ref{coro}) Observe that for $y\in X$, conditions \eqref{effo} and \eqref{atko} together imply condition (a) of the Theorem.
\end{proof}

\begin{proof}(Lemma~\ref{lem1}) Since $F \geq 0$ and $\sigma$ is nondecreasing we see that $y^{\prime}(x) \leq y^{\prime}(0)$ for every $x \geq 0$ (by \eqref{nlvs2}). In addition, for $x_2> x_1 >0$,
$$ y^{\prime}(x_2)-y^{\prime}(x_1) = - \int_{x_1}^{x_2}F(t,y(t))\,d\sigma(t) \leq 0,$$ and so $y^{\prime}(x)$ is nonincreasing. It follows that $y^{\prime}(x)\to L$ where the limit $L \leq y^{\prime}(0)$. The possibility that $-\infty \leq L < 0$ is excluded by the assumption that $y(x) > 0$ for all large $x$. Hence $L$ is finite and non-negative. Suppose that $L \neq 0$. Then, for $\varepsilon > 0$ we can choose $X$ so large that $L-\varepsilon < y^{\prime}(x) < L + \varepsilon$, for every $x \geq X$. Integrating this last expression over $[X, x)$ we get the inequality $y(X) + (x-X)(L-\varepsilon) < y(x) < y(X) + (x-X)(L+\varepsilon)$. Dividing the latter by $x$ and letting $x \to \infty$ we get $y(x) \sim Lx$ as $x \to \infty$. On the other hand, if $L=0$, then the same argument gives us $y(x)/x \to 0$ as $x \to \infty$.
\end{proof}

\begin{proof}(Theorem~\ref{th7}) Assume that \eqref{nlvs2} has a bounded eventually positive solution $y(x)$, with $y(x) > 0$, for all $x \geq x_0.$ An application of Lemma~\ref{lem1} gives that $y^{\prime}(x) \to L$ where $L$ is finite (otherwise $y(x)$ cannot remain bounded at infinity). In addition, passing to the limit as $x \to \infty$ in \eqref{nlvs2}, and rearranging terms, we obtain 
\begin{equation}
\label{ypx}
y^{\prime}(x) = L + \int_{x}^{\infty}F(t,y(t))\,d\sigma(t)
\end{equation}
If $L \neq 0$ then Lemma~\ref{lem1} implies that $y(x)\sim Lx$ as $x \to \infty$ which contradicts the boundedness of $y(x)$. Hence $L=0$. For $x_2>x_1>0$, we integrate \eqref{ypx} over $[x_1,x_2)$ to find that $y(x_2)-y(x_1)>0$ (by our assumptions on $F$ and $\sigma$), that is, $y(x)$ is nondecreasing. Since $y(x)$ is bounded by assumption, we get that $y(x) \to c$ for some finite $c > 0$. Integrating \eqref{ypx} over $[x, \infty)$ and rearranging terms we obtain the existence and finiteness of all integrals involved and, in fact, for all $x\geq x_0$ there holds,
\begin{equation}
\label{fub}
y(x) = c - \int_{x}^{\infty}(t-x)\,F(t,y(t))\,d\sigma(t),
\end{equation}
after an application of Fubini's Theorem. Consolidating our results we have that $0 < y(x_0) < y(x) \leq c$, for all $x \geq x_0$. Observe that the integral \eqref{fub} is finite for $x=x_0$. This, along with the hypothesis that $F(x,\cdot)$ is nondecreasing gives us
\begin{equation}
\label{fub1}
\int_{x_0}^{\infty}(t-x_0)\,F(t,y(x_0))\,d\sigma(t) < \infty,
\end{equation}
and this equivalent to the convergence of \eqref{zn1} with $M=y(x_0)$. Note that we can also replace $M$ by any number smaller than $c$. \\ \\
For the sufficiency we assume that \eqref{zn1} holds for some $M>0$. Fix $A > 0$, $A < M$ and choose $x=a$ so large that $$\int_{a}^{\infty}(t-a)\,F(t,M)\,dt \leq A/2.$$ We set up the iterative scheme
\begin{equation}
\label{iter}
y_{n+1}(x) = A - \int_{x}^{\infty}(t-x)\,F(t, y_{n}(t))\, d\sigma(t), \quad x\geq a,
\end{equation}
with $y_0(x)=A$, for each $x\geq a$. Since $F(t,y_0(t))=F(t,A)\leq F(t,M)$ we get that 
\begin{eqnarray*}
y_{1}(x) &\geq & A - \int_{x}^{\infty}(t-x)\,F(t, M)\, d\sigma(t)\\
& \geq & A - \int_{a}^{\infty}(t-a)\,F(t, M)\, d\sigma(t)\\
&\geq & A - A/2 = A/2.
\end{eqnarray*}
Thus, $A/2 \leq y_1(x) \leq A$ for every $x \geq a$. A similar argument shows that if $A/2 \leq y_n(x) \leq A$ for every $x \geq a$, then the same is true of $y_{n+1}(x)$. An induction argument gives us that
\begin{equation}
\label{a2yn}
A/2 \leq y_n(x) \leq A, \quad x\geq a, n\geq 1.
\end{equation}
Next, we show that each $y_n(x)$ is nondecreasing and the family $\{y_n(x)\}_{n=1}^{\infty}$ is equicontinuous on every interval  $[a,b]$. Let $x_2 > x_1 > a$. Since
\begin{equation}
\label{yn1x}
y_{n+1}(x_2)-y_{n+1}(x_1) = \int_{x_1}^{x_2}(t-x_1)\,F(t,y_n(t))\,d\sigma(t) + \int_{x_2}^{\infty}(x_2-x_1)\,F(t,y_n(t))\,d\sigma(t) ,
\end{equation}
and $F\geq 0$, $\sigma$ is nondecreasing, it follows that the right side of \eqref{yn1x} is non-negative, thus for each $n$, the $y_{n}(x)$ are increasing over $[a, \infty)$. Next, estimating the integrals in \eqref{yn1x} using \eqref{a2yn} and the basic estimates on $A$, we get, for $[x_1, x_2] \in [a,b]$,
\begin{eqnarray*}
|y_{n+1}(x_2)-y_{n+1}(x_1) |&\leq&|x_2-x_1|\left \{\int_{x_1}^{x_2}F(t,M)\,d\sigma(t) + \int_{x_2}^{\infty}F(t,M)\,d\sigma(t)\right \}\\
&\leq &|x_2-x_1|\int_{a}^{\infty}F(t,M)\, d\sigma(t).
\end{eqnarray*}
This last integral, being finite on account of \eqref{zn1}, shows that the family is equicontinuous on $[a,b]$ for every $b > a.$ Thus, passing to a subsequence if necessary, we can say that the limit $y(x) = \lim_{n\to \infty}y_n(x)$ exists and is a continuous function on every interval $[a,b]$.\\

Finally, we show that the limit $y(x)$ is a solution of \eqref{nlvs2} on $[a,\infty)$. Let $\varepsilon > 0$. Rearranging terms in \eqref{iter} we can write, for $x \in [a,b]$,
\begin{eqnarray}
| y_{n+1}(x) - A + \int_{x}^{b}(t-x)\,F(t, y_{n}(t))\, d\sigma(t)| &\leq & \int_{b}^{\infty}(t-x)\,F(t, y_{n}(t))\, d\sigma(t)\nonumber \\
&\leq & \int_{b}^{\infty}(t-a)\,F(t, y_{n}(t))\, d\sigma(t) \nonumber\\
&\leq & \int_{b}^{\infty}(t-a)\,F(t, M)\, d\sigma(t) \\ & < &\varepsilon,\label{ii1}
\end{eqnarray}
provided $b$ is sufficiently large (this is possible on account of \eqref{zn1}). We can now pass to the limit as $n \to \infty$ in \eqref{ii1} to find
\begin{equation*}
| y(x) - A + \int_{x}^{b}(t-x)\,F(t, y(t))\, d\sigma(t)| \leq \varepsilon,
\end{equation*}
holds for every sufficiently large $b$, which is equivalent to saying that 
\begin{equation}
\label{ii0}
y(x) = A - \int_{x}^{\infty}(t-x)\,F(t, y(t))\, d\sigma(t),
\end{equation}
Finally, an application of Fubini's theorem shows that the latter is equivalent to
\begin{equation}
\label{ii2}
y(x) = A - \int_{x}^{\infty}\int_{t}^{\infty}\,F(s, y(s))\, d\sigma(s)\,dt.
\end{equation}
Hence, $y$ is locally absolutely continuous and its (right-) derivative is given at every point $x \geq a$ by differentiating \eqref{ii2}, that is,
\begin{equation}
\label{ii3}
y^{\prime}(x) =  \int_{x}^{\infty}\,F(t, y(t))\, d\sigma(t).
\end{equation}
Writing $y^{\prime}(0)=\int_{0}^{\infty}F(t, y(t))\, d\sigma(t),$ (which necessarily exists and is finite because of \eqref{zn1}) we can rewrite \eqref{ii3} in the form
\begin{eqnarray*}
y^{\prime}(x) &=& y^{\prime}(0) + \int_{x}^{\infty}\,F(t, y(t))\, d\sigma(t) - \int_{0}^{\infty}F(t, y(t))\, d\sigma(t)\\
&=& y^{\prime}(0) - \int_{0}^{x}F(t, y(t))\, d\sigma(t),
\end{eqnarray*}
as desired (see \eqref{nlvs2}).
\end{proof}

\begin{proof}(Corollary~\ref{cor001}) First, we note that the function $F(x,y) := yG(x,y^2)$ satisfies all the conditions of the theorem. In addition, $y$ is a solution of \eqref{nlvs2p} if and only if $-y$ is. Furhermore, \eqref{zn18} is equivalent to \eqref{zn1} for appropriate choices of $c,M$ (indeed, \eqref{zn1} implies \eqref{zn18} with $c=M^2$, and \eqref{zn18} implies \eqref{zn1} with $M=\sqrt{c}$).\\

\noindent{S}ince, for a given solution $y$ of \eqref{nlvs2p} its counterpart $-y$ is also a solution, we can assume without loss of generality that this bounded nonoscillatory solution is eventually positive and so proceed, with no other important changes, as in the proof of the necessity in the theorem to arrive at \eqref{zn18}. The sufficiency proceeds along similar lines.
 \end{proof}

\begin{proof}(Theorem~\ref{th8})  The sufficiency follows the proof of the sufficiency of Theorem~\ref{th7}, which is applicable since \eqref{sts} holds for $\varepsilon =0$, as required by the theorem. Since the solution in Theorem~\ref{th7} is asymptotically a positive constant, it is eventually positive.\\ 

\noindent{F}or the necessity we apply the proof of Lemma~\ref{lem1} to find that if $y(x)$ is eventually positive, say for $x \geq a$, then $c:=\lim_{x\to \infty}y^{\prime}(x) \geq 0$. Thus, with right-derivatives,
$$y^{\prime}(x) = c + \int_{x}^{\infty}F(t,y(t))\,d\sigma(t) \geq \int_{x}^{\infty}F(t,y(t))\,d\sigma(t),$$ and since $y$ is non-decreasing for $x\geq a$ (see the proof of Theorem~\ref{th7}), this gives
\begin{eqnarray*}
{y(x)}^{-\varepsilon}\,y^{\prime}(x) &\geq &\int_{x}^{\infty}{y(t)}^{-\varepsilon}{F(t,y(t))}\,d\sigma(t).
\end{eqnarray*}
Now since $y$ is positive and locally absolutely continuous for $x\geq a$ so is the function $y(x)^{1-\varepsilon}$, since $y(x)$ is bounded away from zero on finite intervals. Hence, for $b > a$, writing $M:=y(a)$,
\begin{eqnarray*}
\int_{a}^{b}{y(x)}^{-\varepsilon}\,y^{\prime}(x)\,dx &\geq &\int_{a}^{b}\int_{x}^{\infty}{y(t)}^{-\varepsilon}\,{F(t,y(t))}\,d\sigma(t)\,dx\\
&\geq &\int_{a}^{b}\int_{x}^{\infty}{M}^{-\varepsilon}\,F(t, M)\,d\sigma(t)\,dx. 
\end{eqnarray*}
Since this is valid for any $b > a$ we can let $b\to \infty$ and simplify the right to get
\begin{eqnarray*}
\int_{a}^{\infty}{y(x)}^{-\varepsilon}y^{\prime}(x)\,dx &\geq &{M}^{-\varepsilon}\int_{a}^{\infty}(t-a)\,F(t,M)\,d\sigma(t).\\
\end{eqnarray*}
The left side is finite since $\varepsilon > 1$ and so the right must be finite, that is, so is \eqref{zn1} for this choice of $M$.
\end{proof}

\begin{proof}(Lemma~\ref{lem2}) By assumption, there exists a number $c>0$ such that for every $x\geq a$, say, we have $y(x)\geq cx$ and $y^{\prime}(x) > 0$ (recall that $y^{\prime}(x)$ tends to a limit as $x \to \infty$, cf., Lemma~\ref{lem1}). Applying \eqref{nlvs2} over $[a,x]$ and rearranging terms we get 
\begin{eqnarray*}
y^{\prime}(a) &= &y^{\prime}(x) + \int_{a}^{x}F(t, y(t))\,d\sigma(t), \quad x\geq a\\
&\geq & \int_{a}^{x}F(t, ct)\,d\sigma(t),
\end{eqnarray*}
for every $x \geq a$. Since $x$ is arbitrary, we can pass to the limit as $x \to \infty$ and thus obtain \eqref{tmt} with $M=c>0$. 
\end{proof}

\begin{proof}(Theorem~\ref{th9}) This follows directly from an application of both Lemma~\ref{lem2} and an application of Theorem~\ref{th5} in the special case where $\sigma$ is nondecreasing.
\end{proof}

\begin{proof} (Theorem~\ref{th15})  We need only sketch the details as they are similar to those included above in the proof of Theorem~\ref{th1}. For $y\in X$ where $X=\{ y\in Y: \|y(x)/f(x)\| \leq 2\}$ consider the map on $X$ defined by \eqref{fnon}. Minor changes in the proof of said Theorem show that, indeed, $T$ is a self-map on $X$ (since $\sigma$ is non-decreasing). In addition, $T$ is a contraction on account of \eqref{t000}. Hence the theorem.
\end{proof}

\begin{proof}(Theorem~\ref{th11}) Define a right-continuous step-function by defining its jumps to be at the integers $n$, of magnitude $\sigma(n)-\sigma(n-0) = - b_n$, for $n\geq 0$  and so $\sigma(t) = {\rm constant}$ in the interval $(n, n+1)$, for $n\in \mathbb{N}$. The integral condition (3) in Theorem~\ref{th5} is equivalent to the above condition on the sum above and the solution $y(x)$ of the Volterra-Stieltjes integro-differential equation \eqref{nlvs2} is such that $y(n):=y_n$ satisfies \eqref{2bdef} for each $n$. The conclusion is a consequence of said theorem.
\end{proof}

\begin{proof} (Theorem~\ref{th12}) We omit the proof as it is similar to that introduced in Theorem~\ref{th11}, with the necessary modifications.
\end{proof}

\end{document}